\documentclass[12pt]{article}

\usepackage{geometry}
\usepackage{amsmath}
\usepackage{amssymb}
\usepackage{bbm}

\usepackage{natbib}

\usepackage{booktabs} 

\usepackage{xcolor}

\usepackage[shortlabels]{enumitem}

\usepackage{graphicx}
\graphicspath{ {figures/} }
\usepackage{placeins}

\newtheorem{definition}{Definition}[section]
\newtheorem{lemma}{Lemma}[section]
\newtheorem{proposition}{Proposition}[section]
\newtheorem{theorem}{Theorem}[section]
\newtheorem{remark}{Remark}[section]
\newtheorem{example}{Example}[section]

\makeatletter
\renewcommand\paragraph{%
	\@startsection{paragraph}
	{4}
	{\z@}
	{3.25ex \@plus1ex \@minus.2ex}
	{-1em}
	{\normalfont\normalsize\bfseries\maybe@addperiod}%
}
\newcommand{\maybe@addperiod}[1]{%
	#1\@addpunct{.}%
}
\makeatother

\usepackage[pagebackref]{hyperref}       
\hypersetup{
	colorlinks=true,
	linkcolor=blue,
	filecolor=blue,
	citecolor=blue,      
	urlcolor=blue,
}
\renewcommand*{\backref}[1]{}
\renewcommand*{\backrefalt}[4]{%
	\ifcase #1 Not cited.%
	\or        Cited on page~#2.%
	\else      Cited on pages~#2.%
	\fi}

\usepackage{cleveref}
\crefname{section}{Section}{Sections}
\crefname{figure}{Figure}{Figures}
\crefname{definition}{Definition}{Definitions}
\crefname{lemma}{Lemma}{Lemmas}
\crefname{proposition}{Proposition}{Propositions}
\crefname{theorem}{Theorem}{Theorems}
\crefname{remark}{Remark}{Remarks}
\crefname{appendix}{Appendix}{Appendices}
\crefname{assumption}{Assumption}{Assumptions}
\crefname{example}{Example}{Examples}

\usepackage[parfill]{parskip}

\usepackage[nottoc,numbib]{tocbibind}

\usepackage[font=footnotesize,labelfont=bf, textfont=it]{caption}

\usepackage{crossreftools}
\DeclareMathOperator*{\argmax}{arg\,max}
\DeclareMathOperator*{\argmin}{arg\,min}
\newcommand{\precprec}{\prec\mathrel{\mkern-5mu}\prec}
\newcommand{\succsucc}{\succ\mathrel{\mkern-5mu}\succ}

\title{Scalarisation-based risk concepts for robust multi-objective optimisation}
\author{
	Ben Tu\thanks{Imperial College London, United Kingdom}
	\and Nikolas Kantas\footnotemark[1]
	\and Robert M. Lee\thanks{BASF SE, Germany}
	\and Behrang Shafei\footnotemark[2]
}
\date{}

\begin{document}
\maketitle
\begin{abstract}
	Robust optimisation is a well-established framework for optimising functions in the presence of uncertainty. The inherent goal of this problem is to identify a collection of inputs whose outputs are both desirable for the decision maker, whilst also being robust to the underlying uncertainties in the problem. In this work, we study the multi-objective case of this problem. We identify that the majority of all robust multi-objective algorithms rely on two key operations: robustification and scalarisation. Robustification refers to the strategy that is used to account for the uncertainty in the problem. Scalarisation refers to the procedure that is used to encode the relative importance of each objective to a scalar-valued reward. As these operations are not necessarily commutative, the order that they are performed in has an impact on the resulting solutions that are identified and the final decisions that are made. The purpose of this work is to give a thorough exposition on the effects of these different orderings and in particular highlight when one should opt for one ordering over the other. As part of our analysis, we showcase how many existing risk concepts can be integrated into the specification and solution of a robust multi-objective optimisation problem. Besides this, we also demonstrate how one can principally define the notion of a robust Pareto front and a robust performance metric based on our ``robustify and scalarise'' methodology. To illustrate the efficacy of these new ideas, we present two insightful case studies which are based on real-world data sets.
\end{abstract}
\section{Introduction}
Robust multi-objective optimisation is concerned with the problem of optimising a vector-valued function under uncertainty. In this paper, we will consider optimising a bounded objective function $f: \mathbb{X} \times \Xi \rightarrow \mathbb{R}^M$, which is defined over a $D$-dimensional input space $\mathbb{X} \subseteq \mathbb{R}^D$ and a $W$-dimensional uncertain parameter space $\Xi \subseteq \mathbb{R}^W$. Here, the inputs $\mathbf{x} \in \mathbb{X}$ are used to denote all the parameters that can be controlled by the decision maker. The uncertain parameters $\boldsymbol{\xi} \in \Xi$ are used to denote all the other variables in the problem, which affect the value of the objective function, but are subject to uncertainty or randomness. Notably, these uncertain parameters might not be controllable. As a result, the goal of robust multi-objective optimisation is to identify a collection of controllable inputs whose distribution of output vectors are in some sense robust to the potential variation in these uncertain parameters. That is, we aim to solve the problem
\begin{equation}
	\underset{{\mathbf{x} \in \mathbb{X}}}{\text{robustmax}} \hspace{1mm} f(\mathbf{x}, \cdot)
	\label{eqn:robust_moo}
\end{equation}
where the `robustmax' operator is some sensible robust generalisation of the vector-valued maximum which has to be defined by the decision maker. 

The majority of existing work in robust multi-objective optimisation has largely focussed on defining the robust optimum via a first principles Pareto ordering style approach where one defines a partial ordering over the space of inputs---see for instance the work by \cite{ide2014mmor,ide2016os} and \cite{botte2019ejoor} for a review. Conceptually, this partial ordering determines whether an input $\mathbf{x} \in \mathbb{X}$ is considered to be more robust than another input $\mathbf{x}' \in \mathbb{X}$ with respect to the objective function $f$. The solution to the corresponding robust optimisation problem is then given by the collection of inputs which are not robustly dominated by any other feasible input. Consequently, one is then typically interested in coming up with a principled computational strategy in order to identify these robust points. The most prevalent approach to accomplish this is the scalarisation strategy. In this approach, one recasts the corresponding robust multi-objective optimisation problem into a sensible collection of scalar-valued optimisation problems; these problems can then be solved using traditional techniques in order to obtain an approximation to the robust set of interest. Although theoretically motivated, the first principles approach is rather rigid in practice because whenever we want to use a new partial ordering, we would have to manually identify and derive a new set of results. In this work, we take a more efficient approach to robust multi-objective optimisation and define robustness directly using scalarisation functions. That is, we define the solution of the robust multi-objective optimisation \eqref{eqn:robust_moo} as being equivalent to solving a collection of general robust scalarised problems that can be readily solved using standard scalar-valued optimisation techniques.

There are the two essential ingredients to our robust multi-objective optimisation strategy: scalarisation and robustification. Scalarisation refers to the procedure of mapping a vector-valued objective function to a family of scalar-valued rewards. Meanwhile, robustification is concerned with the process of quantifying the robust value of a set of possible points. In this work, we focus on the use of risk functionals, otherwise known as risk measures, in order to perform the robustification procedure. These risk-inspired concepts of robustness are very well established in the robust single-objective optimisation literature \citep{quaranta2008job&f,ben-tal2009,bertsimas2011sra,rahimian2019a}, but have not been adequately studied in the multi-objective setting. Our work aims to address this gap by providing a generic methodology for general risk functionals. Similarly, for the scalarisation procedure, our methodology is general and can accommodate for any arbitrary family of scalarisation functions such as the Gerstewitz functional \citep{gerstewitz1983wztl,tammer2020}. Nevertheless, as a working example, we will often use the length scalarisation functions. These are a variant of the Chebyshev scalarisation functions that are known to exhibit many intuitive and useful properties which will be discussed throughout \citep{ishibuchi20092icec,deb20122icec,tu2024aa}. 
 
As identified in preceding papers \citep{fliege2014ejoor, groetzner2022ejoor}, the robustification and scalarisation operations are not commutative in general: the order matters. A core part of our work is a careful analysis of when one should robustify then scalarise (RTS) or whether one should scalarise then robustify (STR). We will highlight later in \cref{sec:versus}, that one should take an RTS approach in the many-time aggregating scenario where one is interested in the aggregated performance at an input over many repeated runs. In contrast, one should take an STR approach in the one-time or non-aggregating scenario where one wishes to obtain consistently strong performance on every individual run.

\subsection{Main contributions} As part of this work, we have made several novel and insightful contributions to the study of robust multi-objective optimisation. Firstly, we have introduced and formalised two general computationally efficient approaches to robust multi-objective optimisation: the RTS approach (\cref{sec:robustify_then_scalarise}) and the STR approach (\cref{sec:scalarise_then_robustify}). As we motivate in \cref{sec:set_theoretic}, these two approaches are very flexible and they encompass many existing robust multi-objective solution strategies as special cases. Secondly, we present a detailed discussion on when one should use an RTS or STR approach (\cref{sec:versus}) and when these approaches are equivalent (\cref{sec:commute}). Thirdly, inspired by the recent results by \cite{tu2024aa}, we showcase in \cref{sec:robust_pareto_fronts} how it is possible to define a valid Pareto front surface which reflects the robust trade-off that is happening under both the RTS and STR paradigm for the length scalarised problems. Specifically, we derive the novel concept of the RTS front (\cref{def:rts_front}) and STR front (\cref{def:str_front}), which are polar surfaces described by the robust length scalarised values. Fourthly, in \cref{sec:performance_metrics}, we present a general class of performance metrics that can be used under the RTS and STR setting in order to assess the performance of any robust multi-objective algorithm. In particular, we generalise the construction of the R2 utility \citep{hansen1998trimm} to the RTS and STR setting. 

\subsection{Structure of the paper} The rest of the paper is organised as follows: In \cref{sec:preliminaries}, we recall the basics of multi-objective optimisation (\cref{sec:moo}), the scalarisation strategy (\cref{sec:scalarisation}) and risk functionals (\cref{sec:risk_functionals}). We then move on to \cref{sec:robust}, where we first introduce the RTS approach (\cref{sec:robustify_then_scalarise}) and STR approach (\cref{sec:scalarise_then_robustify}) and then discuss their relation with existing work (\cref{sec:discussion}). Afterwards, we move on to \cref{sec:robust_pareto_fronts}, where we introduce two generalisations of the Pareto front surface, namely the RTS front and STR front. Besides this, in \cref{sec:performance_metrics}, we introduce the robust R2 utilities, which can be used to quantify the performance of any approximate robust Pareto set. Consequently, in \cref{sec:experiments}, we present some numerical experiments which showcases how our methodology can be used in practice on some real-world data sets. Finally, in \cref{sec:conclusion}, we conclude the paper with a short summary and discussion. All of the proofs for the main results in this paper are presented in \cref{app:proofs}, whilst the scripts that are needed to reproduce the figures and experiments are available in our Github repository: \href{https://github.com/benmltu/scalarize}{\texttt{https://github.com/benmltu/scalarize}}.
\section{Preliminaries}
\label{sec:preliminaries}

\subsection{Multi-objective optimisation}
\label{sec:moo}
Given a bounded vector-valued objective function $g: \mathbb{X} \rightarrow \mathbb{R}^M$, defined over some $D$-dimensional input space $\mathbb{X} \subseteq \mathbb{R}^D$, the standard multi-objective optimisation problem is
\begin{equation}
	\max_{\mathbf{x} \in \mathbb{X}} g(\mathbf{x}),
	\label{eqn:moo}
\end{equation}
where the maximum here is defined using the Pareto partial ordering relations, which we recall below. Note that throughout, we will use $g(\mathbf{x}) \in \mathbb{R}^M$ to denote a generic vector-valued objective; we will use $f(\mathbf{x},\boldsymbol{\xi}) \in \mathbb{R}^M$ when the objective is also dependent on an uncertainty parameter $\boldsymbol{\xi} \in \Xi$. In addition, we will focus entirely on maximisation problems. This can be done without loss of generality as $\min(A) = - \max(-A)$ for any set of vectors $A \subset \mathbb{R}^M$.
\begin{definition}
	[Pareto domination] The weak, strict and strong Pareto domination are denoted by the binary relations $\succeq, \succ$ and $\succsucc$, respectively. We say a vector $\mathbf{y} \in \mathbb{R}^M$ weakly, strictly or strongly Pareto dominates another vector $\mathbf{y'} \in \mathbb{R}^M$ if
	\begin{align*}
		\mathbf{y} \succeq \mathbf{y'} &\iff \mathbf{y} - \mathbf{y'} \in \mathbb{R}_{\geq 0}^M,
		\\
		\mathbf{y} \succ \mathbf{y'} &\iff \mathbf{y} - \mathbf{y'} \in \mathbb{R}_{\geq 0}^M \setminus \{\mathbf{0}_M\},
		\\
		\mathbf{y} \succsucc \mathbf{y'} &\iff \mathbf{y} - \mathbf{y'} \in \mathbb{R}_{> 0}^M,
	\end{align*}
	respectively, where $\mathbf{0}_M \in \mathbb{R}^M$ denotes the $M$-dimensional vector of zeros.
	\label{def:pareto_partial_ordering}
\end{definition}
\begin{definition}
	[Domination region] The Pareto domination (or dominated) region of a set of vectors $A \subseteq \mathbb{R}^M$ is defined as the collection of vectors which dominates (or is dominated by) this set, that is
	\begin{equation*}
		\mathbb{D}_{\diamond}[A] 
		:= \cup_{\mathbf{a} \in A} \{\mathbf{y} \in \mathbb{R}^M: \mathbf{y} \diamond \mathbf{a}\},
	\end{equation*}
	where $\diamond \in \{\succeq, \succ, \succsucc\}$ (or $\diamond \in \{\preceq, \prec, \precprec\}$) denotes a partial ordering relation. In addition, we denote the complement of any domination region by $\mathbb{D}^{C}_{\diamond}[A] := \mathbb{R}^M \setminus \mathbb{D}_{\diamond}[A]$.
	\label{def:domination_region}
\end{definition}

\begin{definition}
	[Pareto optimality] Given a bounded set of vectors $A \subset \mathbb{R}^M$, a point $\mathbf{a} \in A$ is weakly or strictly Pareto optimal if we cannot find another vector $\mathbf{a}' \in A$ which strongly or strictly Pareto dominates it, respectively. The collection of all of the weakly or strictly Pareto optimal points in this set is referred to as the weak or strict Pareto front, $\mathcal{Y}^{\textnormal{weak}}[A] := A \cap \mathbb{D}^{C}_{\precprec}[A]$ and $\mathcal{Y}^{\textnormal{strict}}[A] := A \cap \mathbb{D}^{C}_{\prec}[A]$, respectively.
\end{definition}
Solving the standard multi-objective optimisation problem \eqref{eqn:moo} is equivalent to finding the corresponding optimal set of inputs which maps to the Pareto front of interest. This set of optimal inputs is commonly referred to as the Pareto set. In other words, the standard goal of interest in vector-valued optimisation is to target either the weak Pareto set or the strict Pareto set: $\mathcal{X}_g^{\textnormal{weak}} := g^{-1}(\mathcal{Y}^{\textnormal{weak}}[g(\mathbb{X})])$ or $\mathcal{X}_g^{\textnormal{strict}} := g^{-1}(\mathcal{Y}^{\textnormal{strict}}[g(\mathbb{X})])$, respectively, where $g(\mathbb{X}) := \{g(\mathbf{x})\in\mathbb{R}^M: \mathbf{x} \in \mathbb{X}\}$ is the image of the objective function.

\subsection{Scalarisation perspective}
\label{sec:scalarisation}
The standard strategy that is often used to solve \eqref{eqn:moo} is a scalarisation method. In this approach, we recast the vector-valued optimisation problem of interest into a collection of scalar-valued optimisation problems
\begin{equation}
	\{\max_{\mathbf{x} \in \mathbb{X}} s_{\boldsymbol{\theta}}(g(\mathbf{x})): \boldsymbol{\theta} \in \Theta\},
	\label{eqn:soo}
\end{equation}
where $s_{\boldsymbol{\theta}}: \mathbb{R}^M \rightarrow \mathbb{R}$ is a parametric scalarisation function that is indexed by scalarisation parameters $\boldsymbol{\theta} \in \Theta$. Notably, if one is interested in capturing the Pareto front using \eqref{eqn:soo}, then this would require us to define a suitable\footnote{In theory, one might desire a family of scalarisation functions that can target all Pareto optimal solutions. In practice however, one can usually only recover a finite approximation to the optimal set and therefore this requirement of complete coverage might not be practically necessary.} set of scalarisation functions $\{s_{\boldsymbol{\theta}}\}_{\boldsymbol{\theta} \in \Theta}$ which ensures that the resulting solutions are Pareto optimal. For instance, this former property is usually assured by restricting our attention to scalarisation functions that are monotonically increasing \citep[Definition 2.5]{tu2025sr}. In practice however, a decision maker might not necessarily be interested in just targetting the Pareto front. In particular they might be open to the use of non-monotonic scalarisation functions, which can be used to target desirable points that are not necessarily Pareto optimal. For example, they might consider the use of distance-based scalarisation functions in order to target a set of desired specifications.

\begin{example}
	[Lp scalarisation function] If a decision maker could elicit a desired target vector $\boldsymbol{\upsilon} \in \mathbb{R}^M$, then they might consider setting the scalarisation function to be the negative\footnote{The negative sign is included because we want the points closer to the target vectors to have a higher scalarised value.} weighted $L^p$ distance to this point:
	\begin{equation}
		s^{\textnormal{Lp}}_{(\boldsymbol{\upsilon}, \mathbf{w})}(\mathbf{y}) = - \Bigl(\sum_{m=1}^M (w^{(m)}(\upsilon^{(m)} - y^{(m)}))^p \Bigr)^{1/p}
		\label{eqn:lp_scalarisation}
	\end{equation}
	for some weight vector $\mathbf{w} \in \Delta^{M-1} := \{\mathbf{y} \in \mathbb{R}_{\geq 0}^M: ||\mathbf{y}||_{L^1}=1\}$.  Note that the solution to the corresponding scalarised optimisation problem \eqref{eqn:soo} is not necessarily Pareto optimal if the target vector lies in the feasible space. Specifically, the scalarisation function \eqref{eqn:lp_scalarisation} would not be monotonic over the feasible space and it would actually penalise objective vectors $g(\mathbf{x})\in \mathbb{R}^M$ which Pareto dominate the target vector $\boldsymbol{\upsilon}$.
	\label{eg:lp}
\end{example}

\begin{figure}
	\includegraphics[width=1\linewidth]{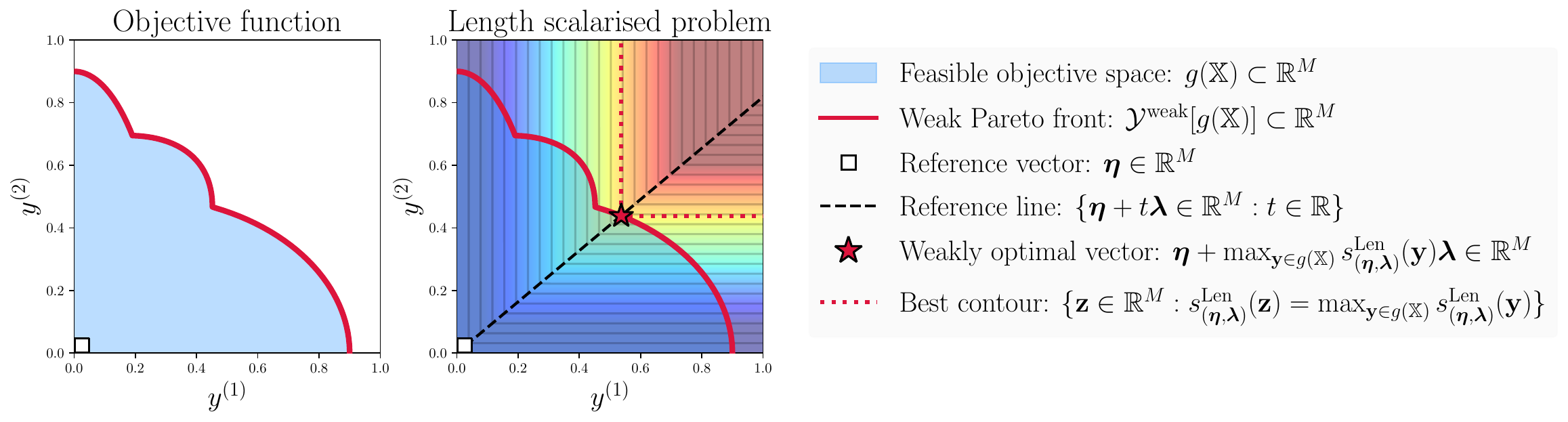}
	\centering
	\caption{An illustration of a length scalarised problem in $M=2$ dimensions.}
	\label{fig:chebyshev_intuition}
\end{figure}

As a working example, we will often use the length scalarisation function, which is a specific reparameterisation of the Chebyshev scalarisation function \citep{ishibuchi20092icec,deb20122icec,tu2024aa}:
\begin{equation}
	s^{\textnormal{Len}}_{(\boldsymbol{\eta}, \boldsymbol{\lambda})}(\mathbf{y})
	= \min_{m=1,\dots,M} \frac{\max (y^{(m)} - \eta^{(m)}, 0)}{\lambda^{(m)}}
	\label{eqn:length_scalarisation}
\end{equation}
for any vector $\mathbf{y} \in \mathbb{R}^M$, where $\boldsymbol{\eta} \in \mathbb{R}^M$ denotes a fixed reference vector and $\boldsymbol{\lambda} \in \mathcal{S}^{M-1}_+$ denotes a positive unit vector that lives in the space of all positive unit vectors $\mathcal{S}^{M-1}_+ := \{\mathbf{y} \in \mathbb{R}_{> 0}^M: ||\mathbf{y}||_{L^2}=1\}$. Conceptually, for any vector $\mathbf{y}$ that strongly dominates the reference vector $\boldsymbol{\eta}$, the length scalarisation function returns the projected length of this vector along the line $\{\boldsymbol{\eta}+t\boldsymbol{\lambda}: t \in \mathbb{R}\}$, that is, $s^{\textnormal{Len}}_{(\boldsymbol{\eta}, \boldsymbol{\lambda})}(\mathbf{y}) = \max\{t > 0: \boldsymbol{\eta} + t \boldsymbol{\lambda} \preceq \mathbf{y}\}$ for any $\mathbf{y} \in \mathbb{D}_{\succsucc}[\{\boldsymbol{\eta}\}]$---see the work by \cite{tu2024aa} for more details for this construction. We now recall a key result from the multi-objective optimisation literature, which states that all weakly Pareto optimal solutions of \eqref{eqn:moo} can be recovered by solving a set of length scalarised problems \eqref{eqn:soo}. For some additional intuition, we have also included \cref{fig:chebyshev_intuition}, which showcases the geometric idea behind this result in two dimensions. 
\begin{theorem}
	[Weak solutions] \cite[Part 2, Theorem 3.4.5]{miettinen1998} Consider a bounded objective function $g: \mathbb{X} \rightarrow \mathbb{R}^M$, and a reference vector $\boldsymbol{\eta} \in \cap_{\mathbf{x} \in \mathbb{X}} \mathbb{D}_{\precprec}[\{g(\mathbf{x})\}]$ which is strongly dominated by the whole feasible objective space. Then, for any weakly Pareto optimal input $\mathbf{z} \in \mathcal{X}_g^{\textnormal{weak}}$, there exists a positive unit vector $\boldsymbol{\lambda} \in \mathcal{S}^{M-1}_+$ such that the input $\mathbf{z}$ is a solution to the scalarised optimisation problem
	\begin{equation}
		\max_{\mathbf{x} \in \mathbb{X}} s^{\textnormal{Len}}_{(\boldsymbol{\eta}, \boldsymbol{\lambda})}(g(\mathbf{x})) 
		= 
		\max_{\mathbf{x} \in \mathbb{X}} \min_{m=1,\dots,M} \frac{g^{(m)}(\mathbf{x}) - \eta^{(m)}}{\lambda^{(m)}}.
		\label{eqn:length_scalarised_problem}
	\end{equation}
	Namely we can choose $\boldsymbol{\lambda} = (g(\mathbf{z}) - \boldsymbol{\eta}) / ||g(\mathbf{z}) - \boldsymbol{\eta}||_{L^2} \in \mathcal{S}^{M-1}_+$.
	\label{thm:chebyshev_scalarisation}
\end{theorem}

\subsection{Risk functionals}
\label{sec:risk_functionals}
\begin{figure}
	\includegraphics[width=1\linewidth]{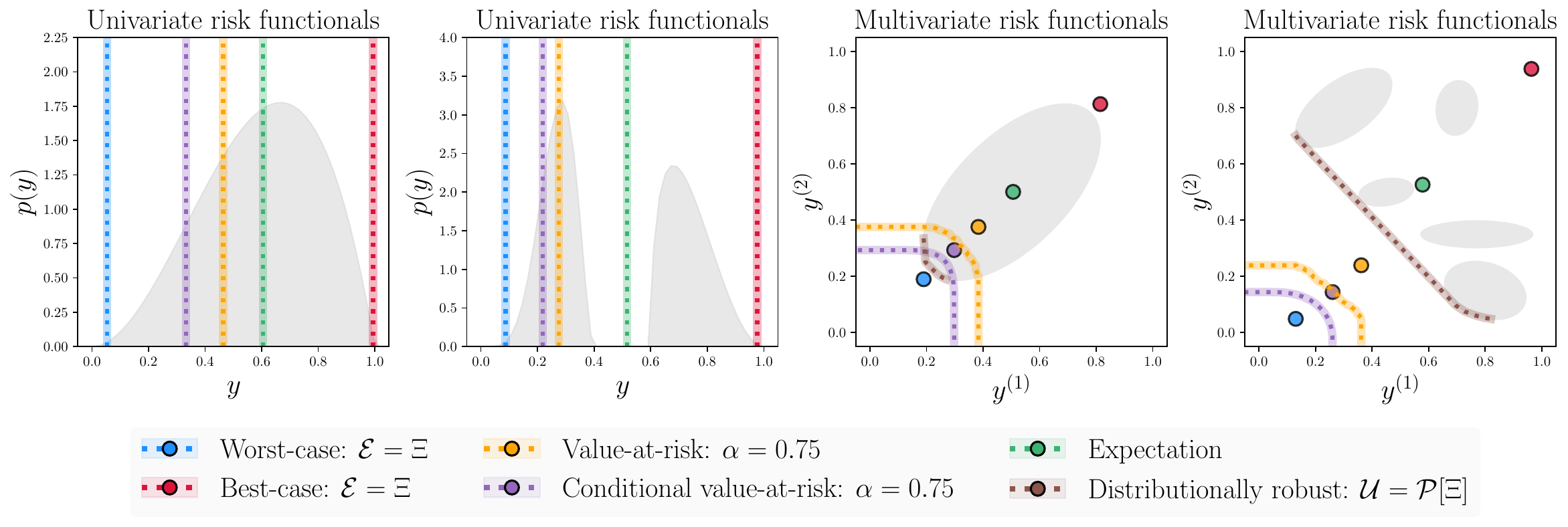}
	\centering
	\caption{An illustration of some univariate and multivariate risk functionals. On the left plots, we consider two univariate probability distributions and plot its associated risk-adjusted values. On the right plots, we consider two uniform distributed random variables and plot its associated risk-adjusted sets. In these multivariate plots, the circles are used to denote the component-wise risk functionals (\cref{eg:multivariate_component_wise}), whilst the lines are used to denote the output of the set-valued risk functionals (\cref{eg:multivariate_var,eg:multivariate_dist_robust,eg:multivariate_pareto_statistics}).}
	\label{fig:risk_functionals}
\end{figure}

In robust optimisation, one generally evaluates the value of an input according to the distribution of its objective values. That is, given an objective function $f: \mathbb{X} \times \Xi \rightarrow \mathbb{R}^M$ the value at any input $\mathbf{x} \in \mathbb{X}$ is typically determined by its corresponding output set
\begin{equation}
	f(\mathbf{x}, \Xi) := \{f(\mathbf{x}, \boldsymbol{\xi}) \in \mathbb{R}^M: \boldsymbol{\xi} \in \Xi\}.
	\label{eqn:output_set}
\end{equation}
Most existing work in robust multi-objective optimisation has focussed on the simple scenario where each uncertain parameter $\boldsymbol{\xi} \in \Xi$ is given equal weight. In this work, we consider a more general treatment where the uncertain parameter is treated as a random variable that is defined sensibly over a standard probability space. More precisely, we suppose that the uncertain parameter $\boldsymbol{\xi} \in \Xi$ is distributed according to some probability density $p \in \mathcal{P}[\Xi]$ that could potentially be known but might in general be unknown. Equipped with this set-up, the value of an input $\mathbf{x} \in \mathbb{X}$ will no longer be determined by its equally weighted output set \eqref{eqn:output_set}, but instead will be defined using a sensible weighted alternative. More specifically, we will consider the use of risk functionals \citep[Part 2]{ruschendorf2013} in order to obtain these proxy set of values. For this reason, we will now give a brief overview on the concepts of a univariate risk functional (\cref{sec:univariate_risk_functionals}) and a multivariate risk functional (\cref{sec:multivariate_risk_functionals}). In \cref{fig:risk_functionals}, we give an illustration of some univariate and multivariate risk functionals for a simple one-dimensional and two-dimensional example, respectively.

\subsubsection{Univariate risk functionals}
\label{sec:univariate_risk_functionals}
A univariate risk functional $\rho$ is any function that maps any real-valued function $h: \Xi \rightarrow \mathbb{R}$ to a real-valued scalar $\rho[h] \in \mathbb{R}$. Intuitively, the value $\rho[h]$ summarises the quality of the output $h(\boldsymbol{\xi}) \in \mathbb{R}$ according to the uncertainty that is inherited from $\boldsymbol{\xi} \in \Xi$. In the following, we present examples of many popular univariate risk functionals that have appeared in the robust optimisation literature.
\begin{example}
	[Extreme cases] Given an uncertainty set $\mathcal{E} \subseteq \Xi$, the worst-case and best-case functional maps any function $h: \Xi \rightarrow \mathbb{R}$ to its corresponding infimum value $\rho^{\textnormal{WorstCase}}[h] := \inf_{\boldsymbol{\xi} \in \mathcal{E}} h(\boldsymbol{\xi})$ and supremum value $\rho^{\textnormal{BestCase}}[h] := \sup_{\boldsymbol{\xi} \in \mathcal{E}} h(\boldsymbol{\xi})$, respectively.
	\label{eg:extreme_cases}
\end{example}
\begin{example}
	[Expectation] Given a probability density $p \in \mathcal{P}[\Xi]$, the expectation functional maps any function $h: \Xi \rightarrow \mathbb{R}$ to its corresponding expectation $\rho^{\textnormal{Exp}}[h] := \mathbb{E}_{p(\boldsymbol{\xi})}[ h(\boldsymbol{\xi})]$.
	\label{eg:expectation}
\end{example}
\begin{example}
	[Value-at-risk] Given a probability density $p \in \mathcal{P}[\Xi]$, the value-at-risk (VaR) functional, at level $\alpha \in (0, 1)$, maps any function $h: \Xi \rightarrow \mathbb{R}$ to its corresponding $(1-\alpha)$-quantile value
	\label{eg:var}
	\begin{equation*}
		\rho^{\textnormal{VaR}_\alpha}[h] := \sup\{y \in \mathbb{R}: \mathbb{E}_{p(\boldsymbol{\xi})}[\mathbbm{1}[h(\boldsymbol{\xi}) \geq y]] \geq \alpha\}.
	\end{equation*}
	\label{eg:value_at_risk}
\end{example}
\begin{example}
	[Distributionally robust] Given an uncertainty set of probability densities $\mathcal{U} \subseteq \mathcal{P}[\Xi]$, the distributionally robust (DR) functional maps any function $h: \Xi \rightarrow \mathbb{R}$ to its distributionally worst-case expected value
	\begin{equation*}
		\rho^{\textnormal{DR}}[h] := \inf_{q \in \mathcal{U}} \mathbb{E}_{q(\boldsymbol{\xi})}[ h(\boldsymbol{\xi})].
	\end{equation*}
	Note that if we set $\mathcal{U} = \mathcal{P}[\mathcal{E}]$, then we recover the worst-case functional $\inf_{q \in \mathcal{P}[\mathcal{E}]} \mathbb{E}_{q(\boldsymbol{\xi})}[ h(\boldsymbol{\xi})] = \inf_{\boldsymbol{\xi} \in \mathcal{E}} h(\boldsymbol{\xi})$. As showcased in existing work \citep{shapiro2017sjo,rahimian2019a}, we can also recover the conditional value-at-risk (CVaR) functional \citep{rockafellar2002job&f} at level $\alpha \in (0, 1)$, with respect to a distribution $p \in \mathcal{P}[\Xi]$ as a special case of the DR functional
	\begin{equation*}
		\rho^{\textnormal{CVaR}_\alpha}[h] := \frac{1}{1-\alpha} \int_{\alpha}^{1}\rho^{\textnormal{VaR}_\tau}[h] d\tau.
	\end{equation*}
	\label{eg:dist_robust}
\end{example}
In the following sections, we will consider the problem of maximising the output of a risk functional. We will often refer to the resulting output of a univariate risk functional as the risk-adjusted or risk-penalised objective value. Naturally, we want these risk functionals to satisfy some sensible properties in order for the corresponding risk-penalised optimisation problem to make sense from a modelling perspective. For this reason, we focus our attention on partially coherent univariate risk functional (\cref{def:partial_coherency}). By design, these functionals satisfy four additional properties which are sufficient to prove many of the results which we state and describe later in \cref{sec:robust_pareto_fronts}. Note that all of the univariate risk functionals listed above (\cref{eg:extreme_cases,eg:expectation,eg:var,eg:dist_robust}) are partially coherent.
\begin{definition}
	[Partial coherency] A univariate risk functional $\rho$ is partially coherent if it satisfies the following four properties:
	\begin{enumerate}[label=(\roman*)]
		\item Normalised: if $h: \Xi \rightarrow \mathbb{R}$ is the zero function, that is $h(\boldsymbol{\xi}) = 0$ for all $\boldsymbol{\xi} \in \Xi$, then $\rho[h] = 0$.
		\item Monotonicity: for any functions $a: \Xi \rightarrow \mathbb{R}$ and $b: \Xi \rightarrow \mathbb{R}$ such that $a(\boldsymbol{\xi}) \geq b(\boldsymbol{\xi})$ for all $\boldsymbol{\xi} \in \Xi$, then $\rho[a] \geq \rho[b]$.
		\item Positively homogeneous: for any positive scalar $\alpha > 0$ and function $h: \Xi \rightarrow \mathbb{R}$, we have that $\rho[\alpha h] = \alpha \rho[h]$.
		\item Translation equivariance: for any constant $c \in \mathbb{R}$ and function $h: \Xi \rightarrow \mathbb{R}$, we have that $\rho[h + c] = \rho[h] + c$.
	\end{enumerate}
	\label{def:partial_coherency}
\end{definition}
\begin{remark}
	[Coherency] A partially coherent univariate risk functional becomes a coherent univariate risk functional, as defined by \cite{artzner1999mf}, if it satisfies an additional property known as sub-additivity: $\rho[a+b] \leq \rho[a] + \rho[b]$ for any functions $a: \Xi \rightarrow \mathbb{R}$ and $b: \Xi \rightarrow \mathbb{R}$. By a well-known duality result \citep[Proposition 4.1]{artzner1999mf}, it is known that any coherent univariate risk functional can be written as a distributionally robust functional (\cref{eg:dist_robust}) with some uncertainty set $\mathcal{U} \subseteq \mathcal{P}[\Xi]$.
\end{remark}
To prove some results later on, we also require that the univariate risk functional are bounded (\cref{def:bounded_univariate_risk}). It is not hard to show that any partially coherent univariate risk functional is indeed bounded. For completeness, we state this result in \cref{lemma:bounded_univariate_risk} and prove it in \cref{app:proofs:lemma:bounded_univariate_risk}.
\begin{definition}
	[Bounded risk] A univariate risk functional $\rho$ is bounded if for any bounded function $h: \Xi \rightarrow \mathbb{R}$, the risk-adjusted value $\rho[h] \in \mathbb{R}$ is bounded.
	\label{def:bounded_univariate_risk}
\end{definition}
\begin{lemma}
	A partially coherent univariate risk functional $\rho$ is bounded.
	\label{lemma:bounded_univariate_risk}
\end{lemma}

\subsubsection{Multivariate risk functionals}
\label{sec:multivariate_risk_functionals}

A multivariate risk functional $\boldsymbol{\rho}$ is any function that maps any vector-valued function $h: \Xi \rightarrow \mathbb{R}^M$ to a closed set of vectors $\boldsymbol{\rho}[h] \subseteq \mathbb{R}^M$. Intuitively, the set of vectors $\boldsymbol{\rho}[h]$ acts as a risk-adjusted summary set of the output set $h(\Xi) \subseteq \mathbb{R}^M$. In the following, we present examples of many popular multivariate risk functionals that have appeared in the robust optimisation literature.

\begin{example}
	[Identity] The identity functional maps any vector-valued function $h: \Xi \rightarrow \mathbb{R}^M$ to its output set $\boldsymbol{\rho}^{\textnormal{Identity}}[h] := \{h(\boldsymbol{\xi}) \in \mathbb{R}^M: \boldsymbol{\xi} \in \Xi\}$.
	\label{eg:multivariate_identity}
\end{example}

\begin{example}
	[Multivariate extreme cases] Given an uncertainty set $\mathcal{E} \subseteq \Xi$, the multivariate worst-case and best-case functional \citep{kuroiwa2012vjom,fliege2014ejoor,ehrgott2014ejoor} maps any vector-valued function $h: \Xi \rightarrow \mathbb{R}^M$ to the singleton set of its corresponding component-wise infimum and supremum value, respectively:
	\begin{align*}
		\boldsymbol{\rho}^{\textnormal{WorstCase}}[h] &:= 
		\{
		(\inf_{\boldsymbol{\xi} \in \mathcal{E}} h^{(1)}(\boldsymbol{\xi}),
		\dots,
		\inf_{\boldsymbol{\xi} \in \mathcal{E}} h^{(M)}(\boldsymbol{\xi}))
		\},
		\\
		\boldsymbol{\rho}^{\textnormal{BestCase}}[h] &:= 
		\{
		(\sup_{\boldsymbol{\xi} \in \mathcal{E}} h^{(1)}(\boldsymbol{\xi}),
		\dots,
		\sup_{\boldsymbol{\xi} \in \mathcal{E}} h^{(M)}(\boldsymbol{\xi}))
		\}.
	\end{align*}
	\label{eg:multivariate_extreme_cases}
\end{example}

\begin{example}
	[Multivariate expectation] Given a probability density $p \in \mathcal{P}[\Xi]$, the multivariate expectation functional maps any vector-valued function $h: \Xi \rightarrow \mathbb{R}^M$ to the singleton set comprised of its corresponding expectation $\boldsymbol{\rho}^{\textnormal{Exp}}[h] := \{\mathbb{E}_{p(\boldsymbol{\xi})}[ h(\boldsymbol{\xi})]\}$. For example, \cite{deb2005emo} considered a special case where the uncertain parameter models any additive input noise arising in \eqref{eqn:moo}, that is $\boldsymbol{\rho}^{\textnormal{Exp}}[f(\mathbf{x}, \cdot)] = \{\mathbb{E}_{p(\boldsymbol{\xi})}[g(\mathbf{x} + \boldsymbol{\xi})]\}$ for any input $\mathbf{x} \in \mathbb{X}$, where $p \in \mathcal{P}[\Xi]$ is the input noise distribution and $g: \mathbb{X} \rightarrow \mathbb{R}^M$ is the vector-valued objective function.
	\label{eg:multivariate_expectation}
\end{example}

\begin{example}
	[Multivariate value-at-risk] Given a probability density $p \in \mathcal{P}[\Xi]$, the multivariate value-at-risk \citep{prekopa2012aor}, at a level $\alpha \in (0, 1)$, maps any vector-valued function $h: \Xi \rightarrow \mathbb{R}^M$ to the set
	\begin{equation*}
		\boldsymbol{\rho}^{\textnormal{VaR}_\alpha}[h] 
		:= \sup\{\mathbf{y} \in \mathbb{R}^M: \mathbb{E}_{p(\boldsymbol{\xi})}[\mathbbm{1}[h(\boldsymbol{\xi}) \succeq \mathbf{y}]] \geq \alpha\},
	\end{equation*}
	where the supremum of a set of vectors is defined with respect to the weak Pareto partial ordering on vectors. At the lowermost case, when $\alpha \downarrow 0$, we obtain the weak Pareto front of the domination region of the essential output set: $\mathcal{Y}^{\textnormal{weak}}[\mathbb{D}_{\preceq}[h(\Xi^{\textnormal{ess}})]]$ where $\Xi^{\textnormal{ess}} := \textnormal{ess supp}(p)$ is the essential support of $p$. In the uppermost case, when $\alpha \uparrow 1$, we obtain the weak Pareto front of the domination region of the worst-case vector: $\mathcal{Y}^{\textnormal{weak}}[\mathbb{D}_{\preceq}[\boldsymbol{\rho}^{\textnormal{WorstCase}}[h]]]$ where $\mathcal{E} = \Xi^{\textnormal{ess}}$.
	\label{eg:multivariate_var}
\end{example}

\begin{example}
	[Multivariate distributionally robust] Given an uncertainty set of probability densities $\mathcal{U} \subseteq \mathcal{P}[\Xi]$, the multivariate distributionally robust functional maps any vector-valued function $h: \Xi \rightarrow \mathbb{R}^M$ to its corresponding distributionally robust worst-case weak Pareto front
	\begin{equation*}
		\boldsymbol{\rho}^{\textnormal{DR}}[h] := \inf_{q \in \mathcal{U}} \mathbb{E}_{q(\boldsymbol{\xi})}[h(\boldsymbol{\xi})],
	\end{equation*}
	where the $\inf$ operator here is defined with respect to the weak Pareto partial ordering on vectors. A special case of this problem has arisen in the work by \citet[Section 5]{bokrantz2017ejoor}, who studied the setting where $\mathcal{U} = \mathcal{P}[\Xi]$. They showed that the corresponding set of vectors is the weak Pareto minimal front of the convex hull of the output set: $-\mathcal{Y}^{\textnormal{weak}}[-\textsc{ConvexHull}[h(\Xi)]]$.
	\label{eg:multivariate_dist_robust}
\end{example}

\begin{example}
	[Component-wise risk] Given a collection of bounded univariate risk functionals $\{\rho^{(1)}$, \dots, $\rho^{(M)}\}$, the component-wise (CW) multivariate risk functional maps any vector-valued function $h: \Xi \rightarrow \mathbb{R}^M$ to the singleton set comprised of its corresponding component-wise risk-adjusted vector
	\begin{equation*}
		\boldsymbol{\rho}^{\textnormal{CW}}[h] := \{(\rho^{(1)}[h^{(1)}], \dots, \rho^{(M)}[h^{(M)}])\},
	\end{equation*}
	where $h(\boldsymbol{\xi}) = (h^{(1)}(\boldsymbol{\xi}), \dots, h^{(M)}(\boldsymbol{\xi})) \in \mathbb{R}^M$ for all $\boldsymbol{\xi} \in \Xi$. The multivariate extreme cases (\cref{eg:multivariate_extreme_cases}) and multivariate expectation (\cref{eg:multivariate_expectation}) are special cases of this component-wise risk functional.
	\label{eg:multivariate_component_wise}
\end{example}

\begin{example}
	[Pareto front statistics] Given a partially coherent univariate risk functional $\rho$ and reference vector $\boldsymbol{\eta} \in \mathbb{R}^M$, the Pareto front surface risk statistic, which we define and derive later in \cref{prop:rho_front_statistics}, maps any vector-valued function $h: \Xi \rightarrow \mathbb{R}^M$ to its corresponding Pareto front risk statistic
	\begin{equation*}
		\boldsymbol{\rho}^{\rho\textnormal{-statistic}}[h]
		:= \rho[\mathcal{Y}^{\textnormal{int}}_{\boldsymbol{\eta}}[\{h(\cdot)\}]]
		= \{
		\boldsymbol{\eta} + 
		\rho[s^{\textnormal{Len}}_{(\boldsymbol{\eta}, \boldsymbol{\lambda})}(h(\cdot)) ]
		\boldsymbol{\lambda}
		\in \mathbb{R}^M: \boldsymbol{\lambda} \in \mathcal{S}_+^{M-1} \}.
	\end{equation*}
	We illustrate an instance of this multivariate risk functional in \cref{fig:risk_functionals} for the CVaR Pareto front surface statistic. Note that this multivariate notion of the CVaR is based on the univariate CVaR (\cref{eg:dist_robust}) and is different from the multivariate generalisation described by \cite{prekopa2012aor} which is much more challenging to define and compute.
	\label{eg:multivariate_pareto_statistics}
\end{example}
For some results later on, we require that the multivariate risk functional of interest is bounded (\cref{def:bounded_multivariate_risk}). Notably, all of the multivariate risk functionals defined above (\cref{eg:multivariate_extreme_cases,eg:multivariate_expectation,eg:multivariate_var,eg:multivariate_dist_robust,eg:multivariate_component_wise,eg:multivariate_pareto_statistics}) are bounded.
\begin{definition}
	[Bounded risk] A multivariate risk functional $\boldsymbol{\rho}$ is bounded if for any bounded vector-valued function $h: \Xi \rightarrow \mathbb{R}^M$, the risk-adjusted set $\boldsymbol{\rho}[h] \subseteq \mathbb{R}^M$ is a bounded set of vectors.
	\label{def:bounded_multivariate_risk}
\end{definition}

\section{Robust multi-objective optimisation}
\label{sec:robust}
Several notions of robustness have already been proposed for the robust multi-objective optimisation problem \eqref{eqn:robust_moo}. Regardless of these different definitions, the most prominent solution strategies that have been proposed to solve these robust problems are all largely based on the scalarisation methodology. We now proceed to introduce the two main solution strategies that will be the primary focus of this work: the robustify then scalarise \eqref{eqn:rts} approach and the scalarise then robustify \eqref{eqn:str} approach. 
\subsection{Robustify then scalarise}
\label{sec:robustify_then_scalarise}
We now introduce a general class of robust multi-objective optimisation problems which we refer to as the robustify then scalarise problems. A scalar-valued optimisation problem belongs to this class if it can be written in the form
\begin{equation}
	\max_{\mathbf{x} \in \mathbb{X}} \max_{\mathbf{y} \in \boldsymbol{\rho}[f(\mathbf{x}, \cdot)]} s_{\boldsymbol{\theta}}(\mathbf{y}),
	\label{eqn:rts}
\end{equation}
where $s_{\boldsymbol{\theta}}: \mathbb{R}^M \rightarrow \mathbb{R}$ denotes any scalarisation function and $\boldsymbol{\rho}$ denotes a multivariate risk functional (\cref{sec:multivariate_risk_functionals}). By design, an RTS problem associates each input $\mathbf{x} \in \mathbb{X}$ with its largest scalarised risk-adjusted value: $\max_{\mathbf{y} \in \boldsymbol{\rho}[f(\mathbf{x}, \cdot)]} s_{\boldsymbol{\theta}}(\mathbf{y}) \in \mathbb{R}$. The solution of an RTS problem is then given by the subset of inputs which map to the largest possible scalarised risk-adjusted value. We refer to these solutions as the RTS robust solutions of $f$ (with respect to $s_{\boldsymbol{\theta}}$ and $\boldsymbol{\rho}$). Following the standard scalarisation methodology (\cref{sec:scalarisation}), the RTS method works by reframing the robust multi-objective optimisation problem \eqref{eqn:robust_moo} as a collection of RTS problems. Intuitively, each of these scalarised problems corresponds to one specific preference $\boldsymbol{\theta} \in \Theta$ and one set of solutions that are RTS robust with respect to this preference. By solving a collection of these problems, one hopes to identify a collection of RTS robust solutions that are desirable for the decision maker. As these RTS problems are all scalar-valued, one can in principle appeal to many standard techniques from single-objective optimisation in order to solve them.

The multivariate risk functional completely characterises the notion of robustness in an RTS problem. At a high level, one can interpret the risk-adjusted output set $\boldsymbol{\rho}[f(\mathbf{x}, \cdot)] \subseteq \mathbb{R}^M$ as being equivalent to the risk-adjusted Pareto front of an input $\mathbf{x} \in \mathbb{X}$. In the special case where the risk-adjusted set only contains one point, the RTS problem reduces down to a standard scalarised optimisation problem \eqref{eqn:soo}, where the vector-valued objective function $g: \mathbb{X} \rightarrow \mathbb{R}^M$ is replaced with some risk-penalised alternative. For instance, this occurs when we use any of the component-wise risk functionals (\cref{eg:multivariate_component_wise}) such as the multivariate expectation (\cref{eg:multivariate_expectation}). Note that we can also recover the original scalarised optimisation problems in \eqref{eqn:soo} by letting $\boldsymbol{\rho}[f(\mathbf{x}, \cdot)] = \{g(\mathbf{x})\}$ for all $\mathbf{x} \in \mathbb{X}$.
\paragraph{Union robust optimisation problem} Recall from \cref{sec:scalarisation} that it is possible to recover any weakly Pareto optimal solution by solving a collection of length scalarised problems---see \cref{thm:chebyshev_scalarisation}.  We will now present the corresponding analogue of this result for the RTS setting. 
\begin{proposition}
	[Union robust weak solutions] Consider a bounded objective function $f: \mathbb{X} \times \Xi \rightarrow \mathbb{R}^M$, a bounded multivariate risk functional $\boldsymbol{\rho}$ and a reference vector $\boldsymbol{\eta} \in \cap_{\mathbf{x} \in \mathbb{X}} \cap_{\mathbf{y} \in \boldsymbol{\rho}[f(\mathbf{x}, \cdot)]} \mathbb{D}_{\precprec}[\{\mathbf{y}\}]$ which is strongly dominated by the whole risk-feasible objective space. Then, for any weakly union robust input $\mathbf{z} \in \mathbb{X}$, there exists a positive unit vector $\boldsymbol{\lambda} \in \mathcal{S}^{M-1}_+$ such that the input $\mathbf{z}$ is a solution to the RTS problem
	\begin{equation}
		\max_{\mathbf{x} \in \mathbb{X}} \max_{\mathbf{y} \in \boldsymbol{\rho}[f(\mathbf{x}, \cdot)]} s^{\textnormal{Len}}_{(\boldsymbol{\eta}, \boldsymbol{\lambda})}(\mathbf{y})
		= 
		\max_{\mathbf{x} \in \mathbb{X}} \max_{\mathbf{y} \in \boldsymbol{\rho}[f(\mathbf{x}, \cdot)]} \min_{m=1,\dots,M} \frac{y^{(m)} - \eta^{(m)}}{\lambda^{(m)}}.
		\label{eqn:rts_length_problem}
	\end{equation}
	\label{prop:robust_chebyshev_scalarisation}
\end{proposition}
To elaborate, \cref{prop:robust_chebyshev_scalarisation} shows that the corresponding RTS length scalarised problem \eqref{eqn:rts_length_problem} can be used in order to target all of the weak solutions to the bounded vector-valued optimisation problem
\begin{equation}
	\max\biggl(\bigcup_{\mathbf{x} \in \mathbb{X}} \boldsymbol{\rho}[f(\mathbf{x}, \cdot)] \biggr),
	\label{eqn:robust_moo_union}
\end{equation}
where $\boldsymbol{\rho}$ is a bounded multivariate risk functional. We refer to this problem as the union robust optimisation problem. Conceptually, this optimisation problem treats every input $\mathbf{x} \in \mathbb{X}$ as being equal to its bounded risk-adjusted output set $\boldsymbol{\rho}[f(\mathbf{x}, \cdot)] \subset \mathbb{R}^M$. An input $\mathbf{x} \in \mathbb{X}$ will then be considered as being weakly or strictly union robust for this problem if one of its output vectors $\mathbf{y} \in \boldsymbol{\rho}[f(\mathbf{x}, \cdot)]$ lies in the overall weak or strict Pareto front of \eqref{eqn:robust_moo_union}, respectively. Through a simple monotonicity argument, one can adapt \cref{thm:chebyshev_scalarisation} in order to show that any weak union robust solution can be obtained as a solution to an RTS problem \eqref{eqn:rts_length_problem}. For completeness, we state the corresponding adaptation of this result in \cref{prop:robust_chebyshev_scalarisation} and prove it in \cref{app:proofs:prop:robust_chebyshev_scalarisation}. 

\subsection{Scalarise then robustify}
\label{sec:scalarise_then_robustify}
We now introduce the other general class of robust multi-objective optimisation problems which we refer to as the scalarise then robustify problems. A scalar-valued optimisation problem belongs to this class if it can be written in the form
\begin{equation}
	\max_{\mathbf{x} \in \mathbb{X}} \rho[s_{\boldsymbol{\theta}}(f(\mathbf{x}, \cdot))],
	\label{eqn:str}
\end{equation}
where $s_{\boldsymbol{\theta}}: \mathbb{R}^M \rightarrow \mathbb{R}$ denotes a scalarisation function and $\rho$ denotes a univariate risk functional (\cref{sec:univariate_risk_functionals}). Conceptually, an STR problem associates each input $\mathbf{x} \in \mathbb{X}$ with its risk-adjusted scalarised value: $\rho[s(f(\mathbf{x}, \cdot))] \in \mathbb{R}$. The solution of an STR problem is then given by the subset of inputs which map to the largest possible risk-adjusted scalarised value. We refer to these solutions as the STR robust solutions of $f$ (with respect to $s_{\boldsymbol{\theta}}$ and $\rho$). Similar with the RTS method (\cref{sec:robustify_then_scalarise}), the STR method works by reframing the robust multi-objective optimisation problem \eqref{eqn:robust_moo} into a collection of STR problems. As before, each scalarised problem corresponds to one specific preference $\boldsymbol{\theta} \in \Theta$ and one set of solutions that are STR robust with respect to this preference. By solving a collection of these problems, one hopes to identify a collection of STR robust solutions that are desirable for the decision maker. Moreover, as these problems are all scalar-valued, one can in theory appeal to many standard techniques from single-objective optimisation in order to solve them.

The univariate risk functional $\rho$ completely characterises the notion of robustness in an STR problem. It determines the risk-adjusted quality of the bounded output set $f(\mathbf{x}, \Xi) \subset \mathbb{R}^M$ according to the preference encoded by the scalarisation function $s_{\boldsymbol{\theta}}$. Note that we can recover the standard scalarised optimisation problems \eqref{eqn:soo} in the special case when we set $\rho[s_{\boldsymbol{\theta}}(f(\mathbf{x}, \cdot))] = s_{\boldsymbol{\theta}}(g(\mathbf{x}))$ for all $\mathbf{x} \in \mathbb{X}$.

\paragraph{STR length scalarised problem} As with the RTS approach \eqref{eqn:rts_length_problem}, one can also define a collection of STR length scalarised problems of the form
	\begin{equation}
	\max_{\mathbf{x} \in \mathbb{X}} \rho[s^{\textnormal{Len}}_{(\boldsymbol{\eta}, \boldsymbol{\lambda})}(f(\mathbf{x}, \cdot))]
	= 
	\max_{\mathbf{x} \in \mathbb{X}} \rho\biggl[\min_{m=1,\dots,M} \frac{f^{(m)}(\mathbf{x}, \cdot) - \eta^{(m)}}{\lambda^{(m)}} \biggr],
	\label{eqn:str_length_problem}
\end{equation}
for any partially coherent univariate risk functional $\rho$, strongly dominated reference vector $\boldsymbol{\eta} \in \cap_{\mathbf{x} \in \mathbb{X}} \cap_{\mathbf{y} \in f(\mathbf{x}, \Xi)} \mathbb{D}_{\precprec}[\{\mathbf{y}\}]$ and positive unit vector $\boldsymbol{\lambda} \in \mathcal{S}^{M-1}_+$. Later on in \cref{sec:robust_pareto_fronts}, we will show that this collection of length scalarised problems is naturally related with the risk statistics associated with the random Pareto front surface of $f$. More precisely, we will see that solving the collection of the STR length scalarised problems \eqref{eqn:str_length_problem} is in some sense equivalent to optimising for a robust Pareto front surface: namely, the STR front (\cref{def:str_front}). In an analogous way, we will also show that solving the collection of RTS length scalarised problems \eqref{eqn:rts_length_problem} is in some sense equivalent to optimising for an RTS front (\cref{def:rts_front}).

\subsection{Discussion}
\label{sec:discussion}
The RTS and STR approach are two general computationally-driven strategies for robust multi-objective optimisation. These two approaches are very flexible and lightweight because they only require the decision maker to specify a sensible family of scalarisation functions and a risk functional in order to define robustness. The former set of functions is determined by the decision maker's internal preferences and specific goals, whilst the latter functional is determined by the uncertainties in the problem and the decision maker's adversity to risk. In this section, we will first discuss the relationships between our methodology and earlier work (\cref{sec:set_theoretic}). We will then look at some natural questions that one might have regarding these two approaches, such as: which approach should one use in practice (\cref{sec:versus}) and when are these two approaches equivalent (\cref{sec:commute})? Moreover, as an aside, we will also relate some of these ideas with the Gerstewitz functional (\cref{sec:gerstewitz}) and existing concepts from multi-objective reinforcement learning (\cref{sec:reinforcement_learning}).

\subsubsection{Related work}
\label{sec:set_theoretic}
\begin{figure}
	\includegraphics[width=1\linewidth]{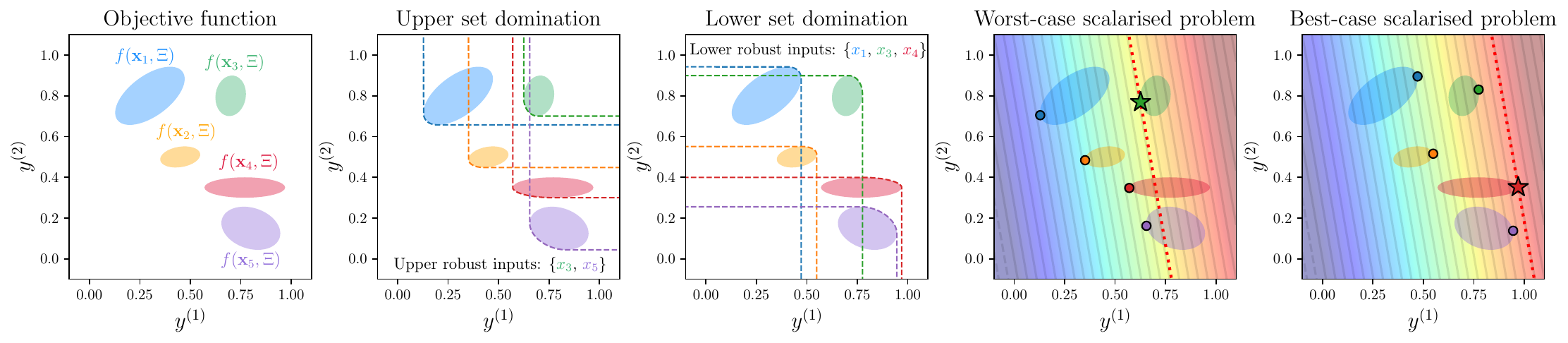}
	\centering
	\caption{An illustration of the upper set and lower set ordered optimisation problem \eqref{eqn:robust_moo_set} based on five output sets. On the left, we highlight the corresponding robust points; on the right, we illustrate how one can solve for these points via the scalarisation-based approaches.}
	\label{fig:set_ordered_robustness}
\end{figure}

Early work on the topic of robust multi-objective optimisation has focussed on treating the problem \eqref{eqn:robust_moo} as being equivalent to the set-valued optimisation problem \citep{khan2015},
\begin{equation}
	\max_{\mathbf{x} \in \mathbb{X}} f(\mathbf{x}, \Xi),
	\label{eqn:robust_moo_set}
\end{equation}
where the maximum here is defined using some set less relations. That is, an input $\mathbf{x} \in \mathbb{X}$ would be considered robust if there does not exist another distinctive input $\mathbf{x}' \in \mathbb{X} \setminus \{\mathbf{x}\}$ whose output set \eqref{eqn:output_set} dominates its own. For example, \cite{ehrgott2014ejoor} proposed using the upper set less relation (\cref{def:set_domination_upper}), whilst \cite{ide2014mmor} and \cite{ide2014fptaa} proposed many other alternatives such as the lower set less relation (\cref{def:set_domination_lower}). As pointed out in these references, all of these proposed approaches reduce down to the standard Pareto partial ordering over vectors (\cref{def:pareto_partial_ordering}) when there is only one possible uncertain parameter $|\Xi|=1$. An illustration of both of these set-orderings is presented on the left of \cref{fig:set_ordered_robustness}.

\begin{definition}
	[Upper set domination] A set $A \subseteq \mathbb{R}^M$ upper dominates a set $B \subseteq \mathbb{R}^M$ if and only if its dominating region is contained in the weak domination region of the other: $A \diamond_u B \iff \mathbb{D}_{\diamond}[A] \subseteq \mathbb{D}_{\succeq}[B]$, for $(\diamond_u, \diamond) \in \{(\succeq_u, \succeq), (\succ_u, \succ), (\succsucc_u, \succsucc)\}$.
	\label{def:set_domination_upper}
\end{definition}

\begin{definition}
	[Lower set domination]	A set $A \subseteq \mathbb{R}^M$ lower dominates a set $B \subseteq \mathbb{R}^M$ if and only if its dominated region is contained in the weak dominated region of the other: $A \diamond_l B \iff \mathbb{D}_{\diamond}[A] \supseteq \mathbb{D}_{\preceq}[B]$, for $(\diamond_l, \diamond) \in \{(\succeq_l, \preceq), (\succ_l, \prec), (\succsucc_l, \precprec)\}$.
	\label{def:set_domination_lower}
\end{definition}
In general, we cannot solve the set-valued optimisation problem \eqref{eqn:robust_moo_set} in its entirety. Instead, most existing optimisation strategies work by identifying some finite approximation to this problem based on some sensible adaptation to the scalarisation approach (\cref{sec:scalarisation}). To elaborate, one usually aims to define and solve a finite collection of scalar-valued optimisation problems whose solutions are robust with respect to the set less relation of interest. Below, we will recall some of these examples for the upper and lower set ordered problem. Moreover, we highlight how these examples are just special cases of our RTS and STR approach.

\paragraph{Objective-wise worst-case problem} One mainstream approach to robust single-objective optimisation is the robust counterpart approach \citep{ben-tal2009}. In this approach, one does not treat the output of an objective function by its entire collection of possibilities \eqref{eqn:output_set}, but instead by some conservative under-estimator of this set: namely, the worst-case value to the output set. The analogue of this approach to the multi-objective setting was presented in earlier work by \cite{kuroiwa2012vjom,fliege2014ejoor} and \cite{ehrgott2014ejoor}. More specifically, they considered the problem of optimising the objective-wise worst-case optimisation problem
\begin{equation}
	\max_{\mathbf{x} \in \mathbb{X}} (\inf_{\boldsymbol{\xi} \in \Xi} f^{(1)}(\mathbf{x}, \boldsymbol{\xi}), \dots, \inf_{\boldsymbol{\xi} \in \Xi} f^{(M)}(\mathbf{x}, \boldsymbol{\xi})).
	\label{eqn:worst_case_problem}
\end{equation}
Evidently, this problem is just an instance of the standard vector-valued optimisation problem \eqref{eqn:moo} and therefore it can be solved using the scalarisation approach described earlier in \cref{sec:scalarisation}. As explained by \cite{ehrgott2014ejoor}, one can obtain solutions to the set-valued optimisation problem \eqref{eqn:robust_moo_set}, under the upper set less relation, by solving different scalarised instances of this worst-case problem \eqref{eqn:worst_case_problem}. Conceptually, this overall computational strategy falls into our class of methods known as the robustify then scalarise approaches. Explicitly, it is the RTS problem where we robustify using the multivariate worst-case functional (\cref{eg:multivariate_extreme_cases}) and scalarise using something like the length scalarisation functions (\cref{thm:chebyshev_scalarisation}).

\paragraph{Extreme case scalarised problems} Alternatively, as suggested in the work by \cite{ehrgott2014ejoor} and \cite{schmidt2019ejoor}, we can also identify some upper set robust inputs to \eqref{eqn:robust_moo_set} by instead solving a collection of worst-case scalarised problems
\begin{equation}
	\max_{\mathbf{x} \in \mathbb{X}} \inf_{\boldsymbol{\xi} \in \Xi} s_{\boldsymbol{\theta}}(f(\mathbf{x}, \boldsymbol{\xi})),
	\label{eqn:str_worst_case}
\end{equation}
for some parametric set of scalarisation functions $s_{\boldsymbol{\theta}}: \mathbb{R}^M \rightarrow \mathbb{R}$. For example, \cite{ehrgott2014ejoor} proposed solving a collection of linearly scalarised problems\footnote{Note that a special instance of these robust linearly scalarised problems also appeared in earlier work by \cite{schottle2009o} and \cite{fliege2014ejoor} for the mean-variance portfolio optimisation problem.} which are defined using the linear scalarisation function $s^{\text{Lin}}_{\mathbf{w}}(\mathbf{y}) = \sum_{m=1}^M w^{(m)} y^{(m)}$ for $\mathbf{y} \in \mathbb{R}^M$, with different weight vectors $\mathbf{w} \in \Delta^{M-1}$. Whereas \cite{schmidt2019ejoor} proposed solving a collection of Chebyshev scalarised problems, which could equivalently be transformed into a collection of length scalarised problems \eqref{eqn:length_scalarisation}. Note that a similar analysis of this kind was also conducted by \cite{ide2014mmor} for the lower set ordered problem. In this latter case, the scalar-valued problems of interest are the best-case scalarised problems
\begin{equation}
	\max_{\mathbf{x} \in \mathbb{X}} \sup_{\boldsymbol{\xi} \in \Xi} s_{\boldsymbol{\theta}}(f(\mathbf{x}, \boldsymbol{\xi})).
	\label{eqn:str_best_case}
\end{equation}
Conceptually, both of these scalarised problems fall into our class of methods known as the scalarise then robustify approaches. Explicitly, the worst-case \eqref{eqn:str_worst_case} and best-case \eqref{eqn:str_best_case} problems scalarises first with some function $s_{\boldsymbol{\theta}}$ and then robustifies using the worst-case and best-case functionals (\cref{eg:extreme_cases}), respectively. On the right of \cref{fig:set_ordered_robustness}, we showcase how both of these approaches can be used along with the linear scalarisation function in order to identify some of the upper and lower set robust solutions to the set-valued optimisation problem \eqref{eqn:robust_moo_set}.

\subsubsection{When to robustify and scalarise?}
\label{sec:versus}
Computationally, the only difference between the RTS and STR approach is the ordering of the operations. Philosophically though, these two approaches differ in both their goals and their interpretations. In the following paragraphs, we elaborate on these differences in some more detail.

\paragraph{Many-time versus one-time philosophy} In general, the RTS approach is suited for problems where the objective function is going to be evaluated many times at the same input and the goal of interest is to maximise the aggregated or summarised performance from these evaluations.  In contrast, the STR approach is suited for problems where the objective function will only be evaluated a finite number of times at any given input and the goal is to maximise the consistent one-time performance at this input. \cref{eg:distance_to_ideal} presents a concrete example of this many-time versus one-time philosophy for the target optimisation problem described in \cref{eg:lp}.
\begin{example}
	[Distance to the ideal vector] Consider the standard Lp scalarisation function \eqref{eqn:lp_scalarisation} with a uniform weight and the expectation functionals (\cref{eg:expectation,eg:multivariate_expectation}). The corresponding RTS problem with this set-up would focus on finding the element that minimises the $L^p$-distance between expected value and the ideal point:
	\begin{equation*}
		\min_{\mathbf{x} \in \mathbb{X}} ||\boldsymbol{\upsilon}- \mathbb{E}_{p(\boldsymbol{\xi})}[f(\mathbf{x}, \boldsymbol{\xi})]||_{L^p}.
	\end{equation*}
	In contrast, the corresponding STR problem would focus on finding the element that minimises the expected $L^p$-distance between the objective function and the ideal point:
	\begin{equation*}
		\min_{\mathbf{x} \in \mathbb{X}} \mathbb{E}_{p(\boldsymbol{\xi})}[||\boldsymbol{\upsilon} - f(\mathbf{x}, \boldsymbol{\xi})||_{L^p}].
	\end{equation*}
	Evidently, these two formulations are different. The former approach is intended for the many-time use setting, where one is interested in identifying an input $\mathbf{x} \in \mathbb{X}$ whose aggregated performance is close to the ideal point. In contrast, the latter approach is designed for the one-time use setting, where one is interested in identifying the input $\mathbf{x} \in \mathbb{X}$ whose sample output vectors are consistently close to the ideal point.
	\label{eg:distance_to_ideal}
\end{example}

\paragraph{Interpretability} In practice, decision makers are interested in making principled decisions that can be justified and explained to other stakeholders. From an interpretation perspective, the RTS approach is appealing because it can easily be interpreted as a standard multi-objective optimisation problem \eqref{eqn:moo} wherein one replaces the objective function with a risk-penalised version \eqref{eqn:robust_moo_union}. The corresponding scalarisation methodology can then be viewed as an operational technique in order to acquire the robust points which are the most desirable. In contrast, the STR approach takes a much more preference-driven strategy that is based largely on the scalarisation methodology. Specifically, one has to first identify a scalarisation function that describes the most preferred total ordering before one can define robustness. This is clearly advantageous in the settings where a preference exists on the outset. But in general, as with the RTS approach, one typically has to keep track of multiple preferences at the beginning before refining them later on. Alternatively, as described earlier in \cref{sec:set_theoretic}, one can also view the STR approach as a computational strategy that is posed in order to solve some set-valued optimisation problem \eqref{eqn:robust_moo_set}.
\subsubsection{Commutative operations}
\label{sec:commute}
There are some special cases where the RTS problem and STR problem are equivalent in the sense that the maximisers are the same, that is
\begin{equation*}
	\argmax_{\mathbf{x} \in \mathbb{X}} \max_{\mathbf{y} \in \boldsymbol{\rho}[f(\mathbf{x}, \cdot)]} s(\mathbf{y})
	= 
	\argmax_{\mathbf{x} \in \mathbb{X}} \rho[s'(f(\mathbf{x}, \cdot))]
\end{equation*}
for some potentially different scalarisation functions $s, s'$ and risk functionals $\rho, \boldsymbol{\rho}$. For example, as described in the previous sections, we can recover the original scalarised optimisation problem \eqref{eqn:soo} as a special case of both the RTS and STR problems. In the following, we present some additional examples where the robustification and scalarisation operation commutes in the above sense. Notably, the primary benefit of this commutativity property is that it gives us a way to get the best of both worlds. That is, we can selectively inherit any desirable features from either the RTS and STR approach in order to justify and solve the robust problem of interest.
\begin{example}
	[Best-case scalarised value] Consider a bounded objective function $f: \mathbb{X} \times \Xi \rightarrow \mathbb{R}^M$ and a scalarisation function $s: \mathbb{R}^M \rightarrow \mathbb{R}$. By simply expanding the definitions, we see that the RTS problem with the identity functional (\cref{eg:multivariate_identity}) and the STR problem with the best-case functional (\cref{eg:extreme_cases}) with $\mathcal{E}=\Xi$ are equivalent:
	\begin{equation*}
		\max_{\mathbf{y} \in f(\mathbf{x}, \Xi)}  s(\mathbf{y})
		= \max_{\boldsymbol{\xi} \in \Xi}s(f(\mathbf{x}, \boldsymbol{\xi}))
	\end{equation*}
	for any input $\mathbf{x} \in \mathbb{X}$.
	\label{eg:equivalence_best_case}
\end{example}

\begin{example}
	[Worst-case length scalarised value] Consider a bounded objective function $f: \mathbb{X} \times \Xi \rightarrow \mathbb{R}^M$, the length scalarisation function \eqref{eqn:length_scalarisation} and the worst-case functionals (\cref{eg:extreme_cases,eg:multivariate_extreme_cases}). By simply expanding the definitions and switching the order of the minimisation and infimum operations, we see that the corresponding RTS and STR problem with this set-up are equivalent:
	\begin{equation*}
		s^{\textnormal{Len}}_{(\boldsymbol{\eta}, \boldsymbol{\lambda})}((\inf_{\boldsymbol{\xi} \in \mathcal{E}} f^{(1)}(\mathbf{x}, \boldsymbol{\xi}), \dots, \inf_{\boldsymbol{\xi} \in \mathcal{E}} f^{(M)}(\mathbf{x}, \boldsymbol{\xi})))
		= \inf_{\boldsymbol{\xi} \in \mathcal{E}} s^{\textnormal{Len}}_{(\boldsymbol{\eta}, \boldsymbol{\lambda})}(f(\mathbf{x}, \boldsymbol{\xi}))
	\end{equation*}
	for any uncertainty set $\mathcal{E} \subseteq \Xi$, reference vector $\boldsymbol{\eta} \in \mathbb{R}^M$, positive unit vector $\boldsymbol{\lambda} \in \mathcal{S}^{M-1}_+$ and input $\mathbf{x} \in \mathbb{X}$.
	\label{eg:equivalence_worst_case}
\end{example}

\begin{example}
	[Linearity of the expectation] Consider a bounded objective function $f: \mathbb{X} \times \Xi \rightarrow \mathbb{R}^M$, the linear scalarisation function and the expectation functionals (\cref{eg:expectation,eg:multivariate_expectation}). By the linearity of the expectation, the corresponding RTS and STR problem with this set-up are equivalent:
	\begin{equation*}
		s^{\textnormal{Lin}}_{\mathbf{w}}(\mathbb{E}_{p(\boldsymbol{\xi})}[f(\mathbf{x}, \boldsymbol{\xi})])
		= \mathbb{E}_{p(\boldsymbol{\xi})}[s^{\textnormal{Lin}}_{\mathbf{w}}(f(\mathbf{x}, \boldsymbol{\xi}))]
	\end{equation*}
	for any distribution $p \in \mathcal{P}[\Xi]$, weight vector $\mathbf{w} \in \Delta^{M-1}$ and input $\mathbf{x} \in \mathbb{X}$.
\end{example}

\begin{example}
	[Value-at-risk of the lengths] Consider a bounded objective function $f: \mathbb{X} \times \Xi \rightarrow \mathbb{R}^M$, the length scalarisation function \eqref{eqn:length_scalarisation} and the VaR functionals (\cref{eg:var,eg:multivariate_var}). By a standard property of the Chebyshev scalarisation function \citep[Theorem 5.1]{daulton2022icml}, the corresponding RTS and STR problem with this set-up are equivalent in the following sense:
	\begin{equation*}
		\max_{\mathbf{y} \in\boldsymbol{\rho}^{\textnormal{VaR}_{\alpha}}[f(\mathbf{x}, \cdot)]}
		s^{\textnormal{Len}}_{(\boldsymbol{\eta}, \boldsymbol{\lambda})}(\mathbf{y})
		=
		\rho^{\textnormal{VaR}_{\alpha}}[s^{\textnormal{Len}}_{(\boldsymbol{\eta}, \boldsymbol{\lambda})}(f(\mathbf{x}, \cdot))]
	\end{equation*}
	for any distribution $p \in \mathcal{P}[\Xi]$, level $\alpha \in (0, 1)$, strongly dominated reference vector $\boldsymbol{\eta} \in \cap_{(\mathbf{x}, \boldsymbol{\xi}) \in \mathbb{X} \times \Xi} \mathbb{D}_{\precprec}[\{f(\mathbf{x}, \boldsymbol{\xi})\}]$,
	positive unit vector $\boldsymbol{\lambda} \in \mathcal{S}^{M-1}_+$ and input $\mathbf{x} \in \mathbb{X}$.
	\label{eg:equivalence_var}
\end{example}
\subsubsection{General nonlinear scalarisation}
\label{sec:gerstewitz}
As showcased in the works by \cite{klamroth2013o, klamroth2017ejoor}, one can unify many different concepts from robust optimisation and stochastic programming via the general class of nonlinear scalarising functionals introduced by Gerstewitz (later Gerth, now Tammer) \citep{gerstewitz1983wztl,tammer2020}. For a topological vector space $V$, the Gerstewitz functional $\varphi_{r-C, k}: V \rightarrow \mathbb{R} \cup \{+\infty\} \cup \{-\infty\}$ is defined by
\begin{equation}
	\varphi_{r-C, k}(v) = \inf\{t \in \mathbb{R}: v \in r - C + tk\},
	\label{eqn:gerstewitz_functional}
\end{equation}
where $C$ is a proper closed subset of $V$, $r \in V$ is a reference point, $k \in V \setminus \{0\}$ is a direction of improvement and $C + k\mathbb{R}_{\geq 0} \subseteq C$. As highlighted in these works, many well-known scalarised minimisation problems are equivalent to minimising an instance of $\varphi_{r-C, k}$ over some subset $W \subseteq V$, that is
\begin{equation}
	\min_{v \in W} \varphi_{r-C, k}(v).
	\label{eqn:gerstewitz_problem}
\end{equation}
For example, in the context of multi-objective optimisation, when $V = \mathbb{R}^M$, we can recover the Chebyshev, $\epsilon$-constraint \citep[Section 4.1]{ehrgott2005} and Pascoletti-Serafini \citep{pascoletti1984jota} scalarised problems as special cases. Similarly, in the context of robust or stochastic single-objective optimisation, when $V = \mathbb{R}^\Xi$, we can recover the worst-case and expectation optimisation problems as special cases. Notably, in the context of both the RTS and STR approaches, we can naturally substitute instances of this functional in place of either the scalarisation function or risk functional (or both).

\paragraph{Length scalarisation function} As the length scalarisation function \eqref{eqn:length_scalarisation} is just a reparameterisation of the Chebyshev scalarisation function, it can also be recovered as a special case of \eqref{eqn:gerstewitz_problem}. More specifically, when $V = \mathbb{R}^M$ we have that
\begin{align*}
	\max_{\mathbf{y} \in Y} s_{(\boldsymbol{\eta}, \boldsymbol{\lambda})}^{\text{Len}}(\mathbf{y})
	&= \max_{\mathbf{y} \in Y} \sup\{t \in \mathbb{R}: \boldsymbol{\eta} + t\boldsymbol{\lambda} \in \mathbf{y} - \mathbb{R}^M_{\geq 0}\}
	\\
	&= - \min_{\mathbf{y} \in Y} \inf\{t \in \mathbb{R}: - \mathbf{y} \in - \boldsymbol{\eta} - \mathbb{R}^M_{\geq 0} + t\boldsymbol{\lambda}\}
	\\
	&= - \min_{\mathbf{y} \in Y} \varphi_{-\boldsymbol{\eta}-\mathbb{R}^M_{\geq 0}, \boldsymbol{\lambda}}(-\mathbf{y})
\end{align*}
where $Y \subseteq \mathbb{D}_{\succsucc}[\{\boldsymbol{\eta}\}]$. It is unclear to us whether the RTS and STR variant of this problem, \eqref{eqn:rts_length_problem} and \eqref{eqn:str_length_problem}, can also be written as special instances of \eqref{eqn:gerstewitz_problem}.

\paragraph{Coherent risk measures} In the context of single-objective optimisation, \cite{quaranta2008job&f} established links between robustification of linear programs and coherent risk measures such as the CVaR. In this work, we only use some of the properties of a coherent risk measure (\cref{def:partial_coherency}) in order to establish the Pareto front risk statistics (\cref{prop:rho_front_statistics}) and the STR front (\cref{prop:str_front}). In general, outside of these specific results, our framework allows one to robustify with risk measures that are not necessarily coherent such as the VaR (\cref{eg:value_at_risk}). 

As explained by \citet[Definition 2.2]{artzner1999mf}, \citet[Section 6]{klamroth2013o} and \citet[Section 15.1.1]{tammer2020}, coherent risk measures (\cref{eg:dist_robust}) can also be written as Gerstewitz functionals \eqref{eqn:gerstewitz_functional}. Generalising our work in terms of these functionals and better understanding its deep connections with coherence and risk are beyond the scope of this paper. These are however very interesting directions for future work. 
\subsubsection{Connections with reinforcement learning}
\label{sec:reinforcement_learning}
Another instance of this RTS versus STR problem has appeared earlier in the topic of multi-objective reinforcement learning \citep{roijers2013j}. In this class of problems, one is interested in identifying policies that maximise the expected cumulative reward associated with a vector-valued Markov decision process. As motivated by \cite{roijers2013j}, scalarisation functions play a prominent role in this problem because they decide which policies one should select and then execute in practice. As highlighted by them, there are two possible scalarised objectives one could define in practice: the scalarised expected return (SER) or the expected scalarised return (ESR). In essence, the SER and ESR approach are just special instances of the RTS and STR approach, respectively. Explicitly, these are the instances where we use the expectation as our risk functional (\cref{eg:expectation,eg:multivariate_expectation}). The philosophical debate on whether one should take an SER or ESR approach is also similar to the RTS versus STR problem (\cref{sec:versus}). Specifically, \citet[Section 8.4]{roijers2013j} says that one should use an SER approach when \textit{``the policy will be used many times and return accumulates across episodes"}. Whereas, an ESR approach should be used when \textit{``the policy will only be used a few times or the return does not accumulate across episodes"}. Note however that in the context of reinforcement learning, the notion of accumulation is that the rewards are summed or averaged. In contrast, in our set-up, we consider the more general notion of aggregation, where the rewards are not necessarily summed but just aggregated into a single or set of summary values via a risk functional.
\section{Robust Pareto fronts}
\label{sec:robust_pareto_fronts}
The Pareto front is an important quantity in standard multi-objective optimisation \eqref{eqn:moo}. It gives us valuable information about the different trade-offs that are happening in the objective space. In practice, this knowledge is typically used by a decision maker in order to help them make an informed and effective decision. So far, much of the work on robust multi-objective optimisation has focussed on classifying and identifying inputs that are deemed as robust in some sense. To the best of our knowledge no work has focussed on identifying and interpreting the robust trade-off that is happening in the objective space. This section addresses this gap in the literature by proposing a principled way to define a robust Pareto front under the RTS and STR perspectives for the length scalarised problems. More precisely, our work leverages ideas from a recent paper by \cite{tu2024aa}, who showed how it is possible to reconstruct an interpolated Pareto front surface from a collection of length scalarised values. This key result was referred to as the polar parameterisation of a Pareto front surface, which we first recall in \cref{sec:polar_parameterisation}. Then we move on to \cref{sec:statistics} where we exploit this parameterisation result in order to define the general concept of a risk statistic associated with a random Pareto front surface (\cref{prop:rho_front_statistics}). Afterwards, in \cref{sec:robust_surfaces}, we present the novel concept of an RTS front (\cref{def:rts_front}) and a STR front (\cref{def:str_front}). Formally, these are the surfaces described by the corresponding RTS and STR length scalarised values, respectively. As we prove in the main results in this section, \cref{prop:rts_front} and \cref{prop:str_front}, both of these surfaces are Pareto front surfaces, under mild assumptions, and therefore are valid robust generalisations to the standard Pareto front trade-off surface. Finally, as an aside, in we present a theoretical result in \cref{sec:bounding_nominal} regarding the ordering of the extreme case Pareto front surfaces.
\subsection{Polar parameterisation}
\label{sec:polar_parameterisation}
To begin with, we recall the definition of a Pareto front surface. Intuitively, this surface is the one obtained by interpolating the strict Pareto front of a set using the weak Pareto partial ordering. As this interpolation could extend indefinitely, we consider truncating this surface at some lower bound defined by a reference vector $\boldsymbol{\eta} \in \mathbb{R}^M$. Afterwards, we recall the definition of a polar surface. 
\begin{definition}
	[Pareto front surface] The Pareto front surface of a bounded set of vectors $A \subset \mathbb{R}^M$, with respect to a reference vector $\boldsymbol{\eta} \in \mathbb{R}^M$, is defined as the truncated weak Pareto front of its weak domination after closure, that is $\mathcal{Y}_{\boldsymbol{\eta}}^{\textnormal{int}}[A] 
	:= \mathcal{Y}^{\textnormal{weak}}[\mathbb{D}_{\preceq}[\textsc{Closure}(A)]] \cap \mathbb{D}_{\succsucc}[\{\boldsymbol{\eta}\}]$. The set of all non-empty Pareto front surfaces with respect to the reference vector $\boldsymbol{\eta} \in \mathbb{R}^M$ is denoted by $\mathbb{Y}_{\boldsymbol{\eta}}^* \subset 2^{\mathbb{R}^M}$.
	\label{def:pareto_front_surface}  
\end{definition}

\begin{definition}
	[Polar surface] A set $A \subseteq \mathbb{R}^M$ is called a polar surface, with respect to a reference vector $\boldsymbol{\eta} \in \mathbb{R}^M$, if and only if there exists a unique non-negative bounded function $r_{\boldsymbol{\eta}, A}: \mathcal{S}^{M-1}_+ \rightarrow \mathbb{R}_{\geq 0}$ such that
	\begin{equation*}
		A = \{\boldsymbol{\eta} + r_{\boldsymbol{\eta}, A}(\boldsymbol{\lambda}) \boldsymbol{\lambda} \in \mathbb{R}^M: \boldsymbol{\lambda} \in \mathcal{S}^{M-1}_+\}.
	\end{equation*}
	The set of all polar surfaces with respect to a reference vector $\boldsymbol{\eta} \in \mathbb{R}^M$ is denoted by $\mathbb{L}_{\boldsymbol{\eta}} \subset 2^{\mathbb{R}^M}$. Moreover, for this set $\mathbb{L}_{\boldsymbol{\eta}}$, we define the projected length function $\ell_{\boldsymbol{\eta}, \boldsymbol{\lambda}}: \mathbb{L}_{\boldsymbol{\eta}} \rightarrow \mathbb{R}_{\geq 0}$ to be the function that returns the projected length along any positive direction $\boldsymbol{\lambda} \in \mathcal{S}^{M-1}_+$, that is $\ell_{\boldsymbol{\eta}, \boldsymbol{\lambda}}[A] = r_{\boldsymbol{\eta}, A}(\boldsymbol{\lambda})$ for any $A \in \mathbb{L}_{\boldsymbol{\eta}}$.
	\label{def:polar_surface}  
\end{definition}
Geometrically, a polar surface is any set that can be obtained by stretching the positive space of unit vectors $\mathcal{S}^{M-1}_+$ along the radial directions, before translating the origin to the reference vector $\boldsymbol{\eta} \in \mathbb{R}^M$. Consequently, this implies that a polar surface can be described entirely by its projected length function $\ell_{\boldsymbol{\eta}, \boldsymbol{\lambda}}$, which computes the amount of stretching that occurs along any positive direction $\boldsymbol{\lambda} \in \mathcal{S}^{M-1}_+$. As shown by \cite{tu2024aa}, and recalled in \cref{thm:polar_parameterisation} below, every Pareto front surface is a polar surface. In particular, it is a polar surface whose projected length function can be written in terms of the length scalarisation function \eqref{eqn:length_scalarisation}. On the left of \cref{fig:polar_parameterisation_risk}, we present a pictorial illustration of this parameterisation result for the Pareto front surface described in \cref{fig:chebyshev_intuition}.

\begin{theorem}
	[Polar parameterisation] \citep[Theorem 3.1]{tu2024aa} For any bounded set of vectors $A \subset \mathbb{R}^M$ and reference vector $\boldsymbol{\eta} \in \mathbb{R}^M$, if the corresponding Pareto front surface is non-empty $\mathcal{Y}_{\boldsymbol{\eta}}^{\textnormal{int}}[A] \neq \emptyset$, then it admits the following polar parameterisation:
	\begin{equation}
		\mathcal{Y}_{\boldsymbol{\eta}}^{\textnormal{int}}[A]
		= \biggl\{\boldsymbol{\eta} + \sup_{\mathbf{a} \in A} s^{\textnormal{Len}}_{(\boldsymbol{\eta}, \boldsymbol{\lambda})}(\mathbf{a}) \boldsymbol{\lambda} \in \mathbb{R}^M: \boldsymbol{\lambda} \in \mathcal{S}_+^{M-1} \biggr\}
		\label{eqn:polar_parameterisation}
	\end{equation}
	where $\ell_{\boldsymbol{\eta}, \boldsymbol{\lambda}}[\mathcal{Y}_{\boldsymbol{\eta}}^{\textnormal{int}}[A]] = \sup_{\mathbf{a} \in A} s^{\textnormal{Len}}_{(\boldsymbol{\eta}, \boldsymbol{\lambda})}(\mathbf{a})$ is the projected length of $A$ along the positive direction $\boldsymbol{\lambda} \in \mathcal{S}_+^{M-1}$.
	\label{thm:polar_parameterisation}
\end{theorem}

\begin{remark}
	[Singleton front] When the Pareto front surface is empty, $\mathcal{Y}_{\boldsymbol{\eta}}^{\textnormal{int}}[A] = \emptyset$, then the right hand side of \eqref{eqn:polar_parameterisation} evaluates to the degenerate polar surface $\{\boldsymbol{\eta}\} \in \mathbb{L}_{\boldsymbol{\eta}}$. This event can happen when the reference vector $\boldsymbol{\eta} \in \mathbb{R}^M$ is set too aggressively in such a way as it weakly dominates the entire feasible objective space.
	\label{rem:singleton}
\end{remark}
\subsection{Risk statistics}
\label{sec:statistics}
\begin{figure}
	\includegraphics[width=1\linewidth]{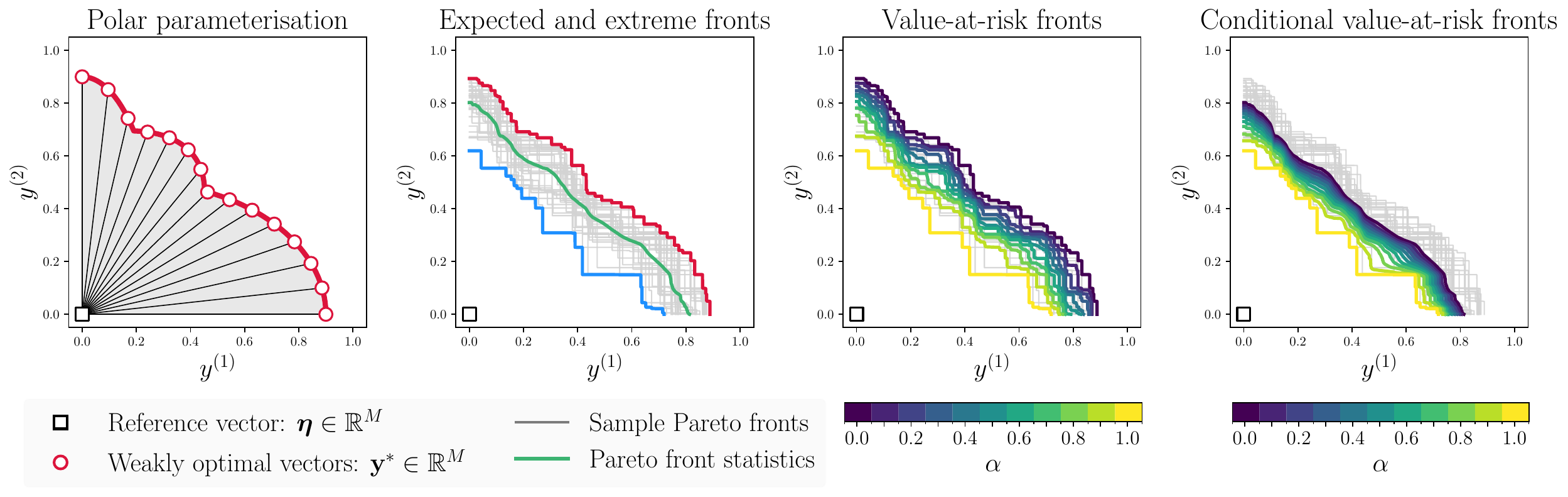}
	\centering
	\caption{We illustrate the polar parameterisation (\cref{thm:polar_parameterisation}) and the Pareto front surface statistics (\cref{prop:rho_front_statistics}) for a simple two-dimensional example. On the left, we plot the polar parameterisation of the Pareto front surface from \cref{fig:chebyshev_intuition}. On the remaining plots, we draw the corresponding risk statistics associated with a finite collection of Pareto front surfaces and the partially coherent univariate risk functionals listed in \cref{sec:univariate_risk_functionals}.}
	\label{fig:polar_parameterisation_risk}
\end{figure}

To apply the polar parameterisation result on the objective function $f: \mathbb{X} \times \Xi \rightarrow \mathbb{R}^M$, with a reference vector $\boldsymbol{\eta} \in \mathbb{R}^M$, we need to ensure that it fulfils the necessary assumptions outlined in \cref{thm:polar_parameterisation}. That is, we need to ensure that the objective function is bounded and that the corresponding Pareto front surface is non-empty. In this work, we only assume that the objective function is bounded. As a result, the polar parameterisation of any bounded set of points could potentially degenerate. Note that this event can however be avoided if the reference vector is set more judiciously (\cref{rem:singleton}). Formally, in this work, we consider the family of polar surfaces
\begin{align}
	\begin{split}
		Y_{\boldsymbol{\eta}, f}^* (\boldsymbol{\xi})
		:= \biggl\{
		\boldsymbol{\eta} + 
		\sup_{\mathbf{x} \in \mathbb{X}} s^{\textnormal{Len}}_{(\boldsymbol{\eta}, \boldsymbol{\lambda})}(f(\mathbf{x}, \boldsymbol{\xi}))
		\boldsymbol{\lambda}
		\in \mathbb{R}^M: \boldsymbol{\lambda} \in \mathcal{S}_+^{M-1} \biggr\},
		\label{eqn:stochastic_pareto_front}
	\end{split}
\end{align}
with $\ell_{\boldsymbol{\eta}, \boldsymbol{\lambda}}[Y_{\boldsymbol{\eta}, f}^* (\boldsymbol{\xi})] = \sup_{\mathbf{x} \in \mathbb{X}} s^{\textnormal{Len}}_{(\boldsymbol{\eta}, \boldsymbol{\lambda})}(f(\mathbf{x}, \boldsymbol{\xi}))$, for any uncertain vector $\boldsymbol{\xi} \in \Xi$ and positive unit vector $\boldsymbol{\lambda} \in \mathcal{S}^{M-1}_+$. If the objective function $f$ is bounded, then this set is either equal to the Pareto front surface $\mathcal{Y}_{\boldsymbol{\eta}}^{\textnormal{int}}[\{f(\mathbf{x}, \boldsymbol{\xi})\}_{\mathbf{x} \in \mathbb{X}}]$, if it is non-empty, or it is the degenerate singleton set $\{\boldsymbol{\eta}\} \in \mathbb{L}_{\boldsymbol{\eta}}$. As shown by \citet[Section 4]{tu2024aa}, one can readily use this polar representation in order to generalise standard statistics from the univariate setting to the space of Pareto front surfaces. For example, one can define the expected Pareto front surface as the polar surface constructed using the expected projected lengths (\cref{eg:expected_pareto_front}).
\begin{example}
	[Expected Pareto front] The expectation of the set \eqref{eqn:stochastic_pareto_front} under some distribution $p \in \mathcal{P}[\Xi]$ is given by the polar surface
	\begin{align*}
		\begin{split}
			\mathbb{E}_{p(\boldsymbol{\xi})}[Y_{\boldsymbol{\eta}, f}^* (\boldsymbol{\xi})]
			:= \{
			\boldsymbol{\eta} + 
			\mathbb{E}_{p(\boldsymbol{\xi})}[\ell_{\boldsymbol{\eta}, \boldsymbol{\lambda}}[Y_{\boldsymbol{\eta}, f}^* (\boldsymbol{\xi})]]
			\boldsymbol{\lambda}
			\in \mathbb{R}^M: \boldsymbol{\lambda} \in \mathcal{S}_+^{M-1} \}.
		\end{split}
	\end{align*}
	As shown by \citet[Section 4.1]{tu2024aa}, if the objective function $f$ is bounded and the Pareto front surfaces $\mathcal{Y}_{\boldsymbol{\eta}}^{\textnormal{int}}[\{f(\mathbf{x}, \boldsymbol{\xi})\}_{\mathbf{x} \in \mathbb{X}}]$ are non-empty almost surely, then this polar surface is a valid Pareto front surface. Otherwise it could degenerate to the singleton set $\{\boldsymbol{\eta}\} \in \mathbb{L}_{\boldsymbol{\eta}}$.
	\label{eg:expected_pareto_front}
\end{example}
We now present a straightforward extension of this result below, which holds for any partially coherent univariate risk functional. The proof of this result is presented in \cref{app:proofs:prop:rho_front_statistics}.
\begin{proposition}
	[Pareto front surface statistics] Consider a bounded objective function $f: \mathbb{X} \times \Xi \rightarrow \mathbb{R}^M$ and a reference vector $\boldsymbol{\eta} \in \mathbb{R}^M$. Then, for any partially coherent univariate risk functional $\rho$, the set
	\begin{align}
		\begin{split}
			\rho[Y_{\boldsymbol{\eta}, f}^*]
			:= \{
			\boldsymbol{\eta} + 
			\rho[\ell_{\boldsymbol{\eta}, \boldsymbol{\lambda}}[Y_{\boldsymbol{\eta}, f}^* (\cdot)] ]
			\boldsymbol{\lambda}
			\in \mathbb{R}^M: \boldsymbol{\lambda} \in \mathcal{S}_+^{M-1} \}
			\label{eqn:rho_front_statistics}
		\end{split}
	\end{align}
	is either a Pareto front surface $\rho[Y_{\boldsymbol{\eta}, f}^*] \in \mathbb{Y}_{\boldsymbol{\eta}}^*$ or the singleton set $\{\boldsymbol{\eta}\} \in \mathbb{L}_{\boldsymbol{\eta}}$.
	\label{prop:rho_front_statistics}
\end{proposition}
On the right of \cref{fig:polar_parameterisation_risk}, we illustrate a number of different risk statistics \eqref{eqn:rho_front_statistics} for a two-dimensional Pareto front surface distribution. In fact, we plot the corresponding empirical estimate for the Pareto surface statistics based on all of the partially coherent univariate risk functionals that were listed in \cref{sec:univariate_risk_functionals}.

\begin{remark}
	[Optimising a risk statistic] Similar to the RTS and STR approach, one could also define the robust optimisation problem \eqref{eqn:robust_moo} as being equal to the problem of optimising a Pareto front surface statistic \eqref{eqn:rho_front_statistics}. That is, we aim to solve a collection of set optimisation problems of the form
	\begin{equation*}
		\max_{X \subseteq \mathbb{X}, |X|\leq P} \rho\biggl[\max_{\mathbf{x} \in X} s^{\textnormal{Len}}_{(\boldsymbol{\eta}, \boldsymbol{\lambda})}(f(\mathbf{x}, \cdot))\biggr],
	\end{equation*}
	for some positive directions $\boldsymbol{\lambda} \in \mathcal{S}^{M-1}_+$, where $P>0$ denotes a cardinality constraint. Individually, each of these problems tries to identify a collection of inputs $X \subseteq \mathbb{X}$, whose collective output, under the preference direction $\boldsymbol{\lambda} \in \mathcal{S}^{M-1}_+$, is of high quality. This problem reduces down to the STR length scalarised problem \eqref{eqn:str_length_problem} when $P=1$. Notably, this is a monotone and submodular set optimisation problem when $P>1$ and therefore is, in general, much more computationally challenging to solve than the corresponding RTS or STR problems---see \cite{tu2025sr} for more details.
	\label{rem:optimise_statistics}
\end{remark}

\subsection{Robust Pareto front surfaces}
\label{sec:robust_surfaces}
The polar parameterisation result in \cref{thm:polar_parameterisation} tells us that we can reconstruct a Pareto front surface from its corresponding collection of maximum length scalarised values. We will now take advantage of this result in order to define the corresponding polar surfaces for the RTS and STR length scalarised problems. More specifically, we define the RTS front and STR front as the polar surface obtained using the RTS length values \eqref{eqn:rts_length_problem} and the STR length values \eqref{eqn:str_length_problem}, respectively.
\begin{definition}
	[RTS front] The RTS front associated with a subset of inputs $X \subseteq \mathbb{X}$, a reference vector $\boldsymbol{\eta} \in \mathbb{R}^M$, a bounded objective function $f: \mathbb{X} \times \Xi \rightarrow \mathbb{R}^M$ and a bounded multivariate risk functional $\boldsymbol{\rho}$ is defined by the set
	\begin{align}
		\begin{split}
			\mathcal{F}_{\boldsymbol{\eta}, f, \boldsymbol{\rho}}^{\textnormal{RTS}}[X]
			:= \biggl\{
			\boldsymbol{\eta} + 
			\sup_{\mathbf{x} \in X} \max_{\mathbf{y} \in \boldsymbol{\rho}[f(\mathbf{x}, \cdot)]} s^{\textnormal{Len}}_{(\boldsymbol{\eta}, \boldsymbol{\lambda})}(\mathbf{y})
			\boldsymbol{\lambda}
			\in \mathbb{R}^M: \boldsymbol{\lambda} \in \mathcal{S}_+^{M-1} \biggr\}.
			\label{eqn:rts_front}
		\end{split}
	\end{align}
	\label{def:rts_front}
\end{definition}
\begin{definition}
	[STR front] The STR front associated with a subset of inputs $X \subseteq \mathbb{X}$, a reference vector $\boldsymbol{\eta} \in \mathbb{R}^M$, a bounded objective function $f: \mathbb{X} \times \Xi \rightarrow \mathbb{R}^M$ and a partially coherent univariate risk functional $\rho$ is defined by the set
	\begin{align}
		\begin{split}
			\mathcal{F}_{\boldsymbol{\eta}, f, \rho}^{\textnormal{STR}}[X]
			:= \biggl\{
			\boldsymbol{\eta} + 
			\sup_{\mathbf{x} \in X} \rho[s^{\textnormal{Len}}_{(\boldsymbol{\eta}, \boldsymbol{\lambda})}(f(\mathbf{x}, \cdot)) ]
			\boldsymbol{\lambda}
			\in \mathbb{R}^M: \boldsymbol{\lambda} \in \mathcal{S}_+^{M-1} \biggr\}.
			\label{eqn:str_front}
		\end{split}
	\end{align}
	\label{def:str_front}
\end{definition}
In \cref{prop:rts_front} and \cref{prop:str_front}, we present the main results in this section which states that these polar surfaces \eqref{eqn:rts_front} and \eqref{eqn:str_front} are indeed valid Pareto front surfaces. Also, we give an explicit representation on what these Pareto front surfaces actually are. In words, the RTS front is the Pareto front surface associated with the union robust problem \eqref{eqn:robust_moo_union}, whilst the STR front is the Pareto front surface associated with the risk statistic of each input's Pareto front surface distribution. The proof of these result are presented in \cref{app:proofs:prop:rts_front} and \cref{app:proofs:prop:str_front}, respectively.

\begin{proposition}
	[RTS front] Consider a bounded objective function $f: \mathbb{X} \times \Xi \rightarrow \mathbb{R}^M$ and reference vector $\boldsymbol{\eta} \in \mathbb{R}^M$. Then, for any bounded multivariate risk functional $\boldsymbol{\rho}$ and any subset of inputs $X \subseteq \mathbb{X}$, the corresponding RTS front is either a Pareto front surface
	\begin{align*}
		\mathcal{F}_{\boldsymbol{\eta}, f, \boldsymbol{\rho}}^{\textnormal{RTS}}[X] 
		&= \mathcal{Y}_{\boldsymbol{\eta}}^{\textnormal{int}}[
		\cup_{\mathbf{x} \in X} \boldsymbol{\rho}[f(\mathbf{x}, \cdot)]]
		= \mathcal{Y}_{\boldsymbol{\eta}}^{\textnormal{int}}[
		\cup_{\mathbf{x} \in X} \mathcal{Y}_{\boldsymbol{\eta}}^{\textnormal{int}}[\boldsymbol{\rho}[f(\mathbf{x}, \cdot)]]]
		\in \mathbb{Y}_{\boldsymbol{\eta}}^*,
	\end{align*}
	or the singleton set $\mathcal{F}_{\boldsymbol{\eta}, f, \boldsymbol{\rho}}^{\textnormal{RTS}}[X] = \{\boldsymbol{\eta}\} \in \mathbb{L}_{\boldsymbol{\eta}}$. 
	\label{prop:rts_front}
\end{proposition}

\begin{figure}
	\includegraphics[width=1\linewidth]{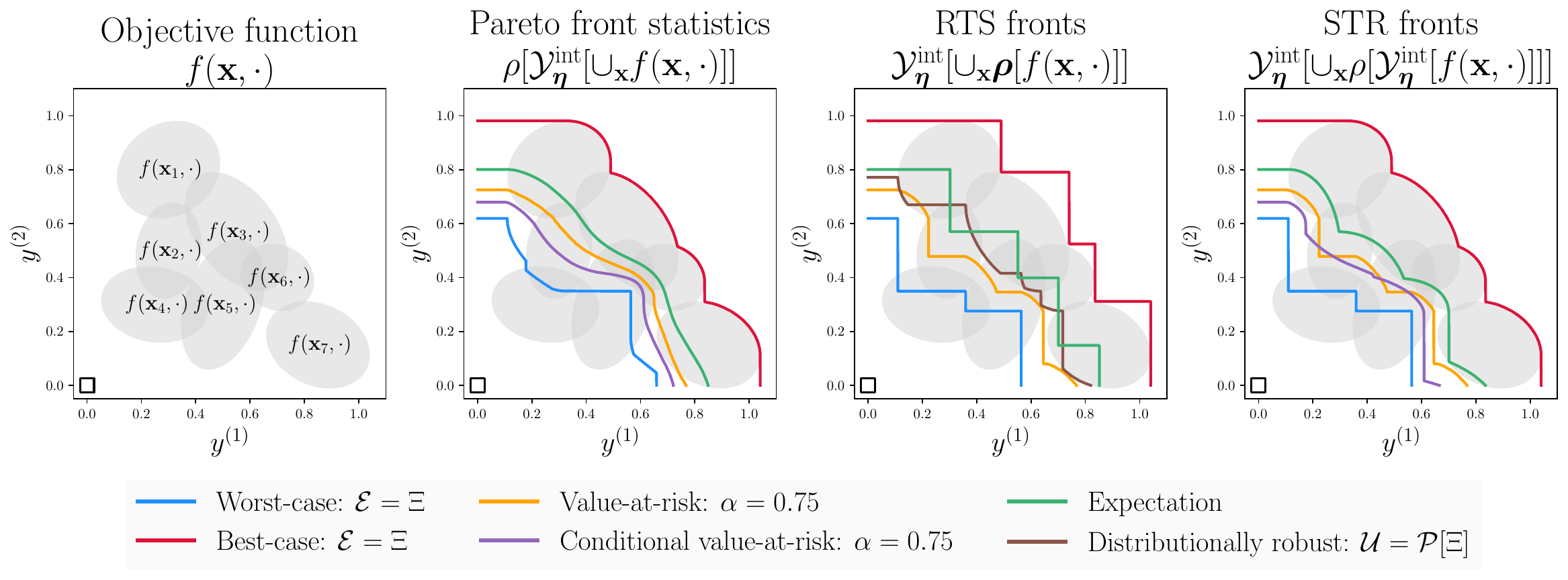}
	\centering
	\caption{An illustration of the different Pareto front surfaces that can be computed from an $M=2$ dimensional objective function. On the left, we plot the output sets $f(\mathbf{x}, \Xi) \subset \mathbb{R}^M$ for each input $\mathbf{x} \in \mathbb{X}$, which are all assumed to be uniformly distributed in some ellipsoid. On the remaining plots, we draw the Pareto front surfaces associated with the Pareto front statistics, RTS fronts and STR fronts for the different risk functionals described in \cref{sec:risk_functionals}.}
	\label{fig:robust_fronts}
\end{figure}

\begin{proposition}
	[STR front] Consider a bounded objective function $f: \mathbb{X} \times \Xi \rightarrow \mathbb{R}^M$ and reference vector $\boldsymbol{\eta} \in \mathbb{R}^M$. Then, for any partially coherent univariate risk functional $\rho$ and any subset of inputs $X \subseteq \mathbb{X}$, the corresponding STR front is either a Pareto front surface
	\begin{align*}
		\mathcal{F}_{\boldsymbol{\eta}, f, \rho}^{\textnormal{STR}}[X] 
		= \mathcal{Y}_{\boldsymbol{\eta}}^{\textnormal{int}}[\cup_{\mathbf{x} \in X} \rho[\mathcal{Y}_{\boldsymbol{\eta}}^{\textnormal{int}}[f(\mathbf{x}, \cdot)]] \in \mathbb{Y}_{\boldsymbol{\eta}}^*,
	\end{align*}
	or the singleton set $\mathcal{F}_{\boldsymbol{\eta}, f, \rho}^{\textnormal{STR}}[X] = \{\boldsymbol{\eta}\} \in \mathbb{L}_{\boldsymbol{\eta}}$.
	\label{prop:str_front}
\end{proposition}
In \cref{fig:robust_fronts}, we illustrate these risk-adjusted Pareto front surfaces for a simple two-dimensional example with the risk functionals introduced in \cref{sec:risk_functionals}. We also compare these Pareto front surfaces with the risk statistics described in \cref{prop:rho_front_statistics}. Overall, we see a great variety between the Pareto front surfaces illustrated within each plot; the differences among the plots however are very subtle. Conceptually, the Pareto front surface statistics looks at the aggregated performance at each scalarised value and therefore has a noticeably smoother surface. In contrast, the RTS and STR front looks at the individual performance at each scalarised level and therefore is a bit more bumpy. In particular, the RTS front treats every input $\mathbf{x} \in X$ according to the Pareto front surface of its risk-adjusted output set $\mathcal{Y}_{\boldsymbol{\eta}}^{\textnormal{int}}[\boldsymbol{\rho}[f(\mathbf{x}, \cdot)]]$, whilst the STR front treats every input $\mathbf{x} \in X$ according to a risk statistic of its random Pareto front surface $\rho[\mathcal{Y}_{\boldsymbol{\eta}}^{\textnormal{int}}[f(\mathbf{x}, \cdot)]]$. These two formulations are different in general. However, as highlighted in \cref{sec:commute}, there are special cases when they are equivalent, such as: the worst-case setting (blue lines) and the VaR setting (orange lines).

\subsection{Bounding a nominal surface}
\label{sec:bounding_nominal}
In some practical settings, the uncertain parameter $\boldsymbol{\xi} \in \Xi$ denotes some nuisance parameter that obscures the true value of the objective function. That is, we suppose that the true value of the objective function appears when this parameter realises its nominal value of $\boldsymbol{\xi}^* \in \Xi$. Under this set-up, the traditional goal for a practitioner would be to identify the Pareto front surface at this nominal value: $\mathcal{Y}_{\boldsymbol{\eta}}^{\textnormal{int}}[\{f(\mathbf{x}, \boldsymbol{\xi}^*)\}_{\mathbf{x} \in \mathbb{X}}]$. In general, one cannot do this exactly. Instead, one typically tries to identify a bound for this surface which holds with high probability. In this section, we analyse the gap between the upper and lower bounding surface obtained using the extreme case functionals (\cref{eg:extreme_cases,eg:multivariate_extreme_cases}). Our main result in this section is \cref{prop:extreme_case_bounds}. This result states that the gap between the best-case and worst-case Pareto front surfaces, for the different robust variants, are ordered. The smallest gap is attained by the Pareto front surface statistics, followed by the STR fronts and then the RTS fronts. This result is illustrated in \cref{fig:robust_fronts} using the best-case fronts (red lines) and worst-case fronts (blue lines). We will now set up the notation that is required in order to state this result.

For any bounded objective function $f: \mathbb{X} \times \Xi \rightarrow \mathbb{R}^M$ and strongly dominated reference vector $\boldsymbol{\eta} \in \cap_{\mathbf{x} \in \mathbb{X}} \cap_{\mathbf{y} \in f(\mathbf{x}, \Xi)} \mathbb{D}_{\precprec}[\{\mathbf{y}\}]$, we can write the target Pareto front surface in its polar parameterised form \eqref{eqn:stochastic_pareto_front}. By definition, this implies that the projected lengths of this Pareto front surface $Y_{\boldsymbol{\eta}, f}^*( \boldsymbol{\xi}^*)$, along any positive direction $\boldsymbol{\lambda} \in \mathcal{S}^{M-1}_+$, can be written as
\begin{equation*}
	\ell_{\boldsymbol{\eta}, \boldsymbol{\lambda}}[Y_{\boldsymbol{\eta}, f}^*(\boldsymbol{\xi}^*)] 
	= \sup_{\mathbf{x} \in \mathbb{X}} \min_{m=1,\dots,M} z^{(m)}_{\boldsymbol{\eta}, \boldsymbol{\lambda}}(f^{(m)}(\mathbf{x}, \boldsymbol{\xi}^*)),
\end{equation*}
where the term $z^{(m)}_{\boldsymbol{\eta}, \boldsymbol{\lambda}}(y) = (y - \eta^{(m)}) /\lambda^{(m)}$ is used to denote the $m$-th argument in the length scalarisation function \eqref{eqn:length_scalarisation} for any number $y \in \mathbb{R}$. It can be easily shown that this Pareto front surface can be bounded from above and below by the extreme cases of the Pareto front surface statistics (denoted by PS), RTS fronts and STR fronts. Mathematically, this statement can be formalised by the following three inequalities:
\begin{align*}
	\underbrace{\inf_{\boldsymbol{\xi} \in \Xi} \sup_{\mathbf{x} \in \mathbb{X}} \min_{m=1,\dots,M} z^{(m)}_{\boldsymbol{\eta}, \boldsymbol{\lambda}}(f^{(m)}(\mathbf{x}, \boldsymbol{\xi}))}_{l_{\boldsymbol{\eta}, \boldsymbol{\lambda}, f}^{\textnormal{PS}}}
	\leq \ell_{\boldsymbol{\eta}, \boldsymbol{\lambda}}[Y_{\boldsymbol{\eta}, f}^*(\boldsymbol{\xi}^*)]
	\leq 
	\underbrace{\sup_{\boldsymbol{\xi} \in \Xi} \sup_{\mathbf{x} \in \mathbb{X}} \min_{m=1,\dots,M} z^{(m)}_{\boldsymbol{\eta}, \boldsymbol{\lambda}}(f^{(m)}(\mathbf{x}, \boldsymbol{\xi}))}_{u_{\boldsymbol{\eta}, \boldsymbol{\lambda}, f}^{\textnormal{PS}}},
	\\
	\underbrace{\sup_{\mathbf{x} \in \mathbb{X}} \min_{m=1,\dots,M} z^{(m)}_{\boldsymbol{\eta}, \boldsymbol{\lambda}}(\inf_{\boldsymbol{\xi} \in \Xi} f^{(m)}(\mathbf{x}, \boldsymbol{\xi}))}_{l_{\boldsymbol{\eta}, \boldsymbol{\lambda}, f}^{\textnormal{RTS}}}
	\leq \ell_{\boldsymbol{\eta}, \boldsymbol{\lambda}}[Y_{\boldsymbol{\eta}, f}^*(\boldsymbol{\xi}^*)]
	\leq 
	\underbrace{\sup_{\mathbf{x} \in \mathbb{X}} \min_{m=1,\dots,M} z^{(m)}_{\boldsymbol{\eta}, \boldsymbol{\lambda}}(\sup_{\boldsymbol{\xi} \in \Xi} f^{(m)}(\mathbf{x}, \boldsymbol{\xi}))}_{u_{\boldsymbol{\eta}, \boldsymbol{\lambda}, f}^{\textnormal{RTS}}},
	\\
	\underbrace{\sup_{\mathbf{x} \in \mathbb{X}} \inf_{\boldsymbol{\xi} \in \Xi} \min_{m=1,\dots,M} z^{(m)}_{\boldsymbol{\eta}, \boldsymbol{\lambda}}(f^{(m)}(\mathbf{x}, \boldsymbol{\xi}))}_{l_{\boldsymbol{\eta}, \boldsymbol{\lambda}, f}^{\textnormal{STR}}}
	\leq \ell_{\boldsymbol{\eta}, \boldsymbol{\lambda}}[Y_{\boldsymbol{\eta}, f}^*(\boldsymbol{\xi}^*)]
	\leq 
	\underbrace{\sup_{\mathbf{x} \in \mathbb{X}} \sup_{\boldsymbol{\xi} \in \Xi} \min_{m=1,\dots,M} z^{(m)}_{\boldsymbol{\eta}, \boldsymbol{\lambda}}(f^{(m)}(\mathbf{x}, \boldsymbol{\xi}))}_{u_{\boldsymbol{\eta}, \boldsymbol{\lambda}, f}^{\textnormal{STR}}},
\end{align*}
which all hold simultaneously for every uncertain parameter $\boldsymbol{\xi}^* \in \Xi$. We see clearly that the only difference between these bounds are the ordering of the three main operations: maximisation over the inputs, robustification over the uncertain parameters and scalarisation of the objective.
\begin{proposition}
	[Extreme case bounds] For any bounded objective function $f: \mathbb{X} \times \Xi \rightarrow \mathbb{R}^M$ and strongly dominated reference vector $\boldsymbol{\eta} \in \cap_{\mathbf{x} \in \mathbb{X}} \cap_{\mathbf{y} \in f(\mathbf{x}, \Xi)} \mathbb{D}_{\precprec}[\{\mathbf{y}\}]$, the following inequality holds:
	\begin{equation*}
		0 \leq u_{\boldsymbol{\eta}, \boldsymbol{\lambda}, f}^{\textnormal{PS}} - l_{\boldsymbol{\eta}, \boldsymbol{\lambda}, f}^{\textnormal{PS}}
		\leq u_{\boldsymbol{\eta}, \boldsymbol{\lambda}, f}^{\textnormal{STR}} - l_{\boldsymbol{\eta}, \boldsymbol{\lambda}, f}^{\textnormal{STR}}
		\leq u_{\boldsymbol{\eta}, \boldsymbol{\lambda}, f}^{\textnormal{RTS}} - l_{\boldsymbol{\eta}, \boldsymbol{\lambda}, f}^{\textnormal{RTS}},
	\end{equation*}
	for any positive direction $\boldsymbol{\lambda} \in \mathcal{S}^{M-1}_+$.
	\label{prop:extreme_case_bounds}
\end{proposition}
\cref{prop:extreme_case_bounds} states that the ordering of these three operations can be exploited in order to get a smaller bounding interval. The proof of this result is presented in \cref{app:proofs:prop:extreme_case_bounds}. Note that in the special case when the set of uncertain parameters is a singleton set $\Xi=\{\boldsymbol{\xi}^*\}$, all of these extreme case Pareto front surfaces are equal. In other words, the difference between these three approaches only begins to transpire when the uncertainty set grows. \cref{prop:extreme_case_bounds} gives us an ordering on the range of Pareto front surfaces one can obtain. It does not however give us any information regarding which method is the most effective to identify or estimate this uncertain parameter $\boldsymbol{\xi}^* \in \Xi$. This latter problem turns out to be case specific as it relies heavily on the assumptions that are placed on the objective function and the form of the uncertainty.
\section{Robust performance metrics}
\label{sec:performance_metrics}
Performance metrics are often used in multi-objective optimisation in order to quantify the quality of a Pareto front approximation. As showcased in a recent survey by \cite{audet2021ejoor}, there exists many different types of performance metrics in the multi-objective literature, all of which try to measure different facets of the Pareto front approximation. Very little work has focussed on extending these existing performance metrics to the robust setting \citep{mirjalili2015saeca}. More specifically, in this latter setting, one wishes to quantify the quality of a set of proposed inputs $X \subseteq \mathbb{X}$, based on the corresponding set of output sets $\{f(\mathbf{x}, \Xi)\}_{\mathbf{x} \in X}$. Naturally, if the individual output sets are singletons, $|f(\mathbf{x}, \Xi)| = 1$, then we reduce down to the standard multi-objective set-up. In this section, we will now address this gap in the literature, and present two family of robust multi-objective performance metrics based on the RTS and STR approaches. Notably, we generalise the R2 utilities \citep{hansen1998trimm}, which are a general family of multi-objective performance metrics that exhibit many favourable properties and contains many of the most widely-used performance metrics as special cases \citep{tu2025sr}. As explicit examples, we will proceed to show how one can generalise the well-known inverted generational distance (\cref{sec:robust_igd}) and hypervolume indicator (\cref{sec:robust_hypervolume}) to the robust setting.
\subsection{Robust R2 utilities}
\label{sec:robust_r2}
Formally, a multi-objective performance metric is any utility function $U: 2^{\mathbb{R}^M} \rightarrow \mathbb{R}$, which maps a set of vectors to a scalar. We adopt the convention that larger utility values signify better performance. The R2 utilities are a special class of utility functions that is defined with respect to a family of scalarisation functions $\{s_{\boldsymbol{\theta}}: \mathbb{R}^M \rightarrow \mathbb{R}\}_{\boldsymbol{\theta} \in \Theta}$ and a probability distribution over the corresponding scalarisation parameters $p \in \mathcal{P}[\Theta]$. Specifically, equipped with this set-up, the corresponding R2 utility of a set of vectors $Y \subset \mathbb{R}^M$ is equal to the weighted average maximum scalarised value
\begin{equation}
	U[Y] := \mathbb{E}_{p(\boldsymbol{\theta})}\biggl[
	\sup_{\mathbf{y} \in Y} s_{\boldsymbol{\theta}}(\mathbf{y})
	\biggr].
	\label{eqn:r2_utility}
\end{equation}
Note that this is the general form of an R2 utility which also accounts for the setting where the Pareto front approximation $Y \subset \mathbb{R}^M$ is continuous. Nevertheless, in practice, the Pareto front approximation is typically finite and therefore the supremum operation is often replaced with the maximum. At an intuitive level, an R2 utility just computes a weighted score associated with a collection of scalarised problems \eqref{eqn:soo}. By design, this family of utility functions is very flexible and it gives us any easy way to define and measure different notions of performance in the objective space. In theory, one could take advantage of this flexibility in order to define R2 utilities which adequately cater for the decision maker's internal preferences.

\paragraph{Robust R2 Utility} Inspired by the construction of the R2 utilities \eqref{eqn:r2_utility}, we now define the corresponding RTS and STR extension of the R2 utilities by replacing the maximum scalarised values arising in \eqref{eqn:soo} with the RTS scalarised values \eqref{eqn:rts} and STR scalarised values \eqref{eqn:str}, respectively. More precisely, we define an RTS-R2 utility function $R_{f, \boldsymbol{\rho}}^{\text{RTS}}: 2^{\mathbb{X}} \rightarrow \mathbb{R}$ and a STR-R2 utility function $R_{f, \rho}^{\text{STR}}: 2^{\mathbb{X}} \rightarrow \mathbb{R}$ by the following equations:
\begin{align}
	R^{\text{RTS}}_{f, \boldsymbol{\rho}}[X] 
	&:= \mathbb{E}_{p(\boldsymbol{\theta})}\biggl[
	\sup_{\mathbf{x} \in X} \max_{\mathbf{y} \in \boldsymbol{\rho}[f(\mathbf{x}, \cdot)]} s_{\boldsymbol{\theta}}(\mathbf{y})
	\biggr],
	\label{eqn:rts_utility}
	\\
	R^{\text{STR}}_{f, \rho}[X] &:= \mathbb{E}_{p(\boldsymbol{\theta})}\biggl[
	\sup_{\mathbf{x} \in X} \rho[s_{\boldsymbol{\theta}}(f(\mathbf{x}, \cdot))]
	\biggr],
	\label{eqn:str_utility}
\end{align}
for any set of inputs $X \subseteq \mathbb{X}$. Note that in contrast to the original R2 utilities in \eqref{eqn:r2_utility}, these robust R2 utilities assess the quality of a set of inputs instead of a set of outputs. Formally, these robust performance metrics are defined with respect to an objective function and a risk functional. This construction makes sense from an intuitive perspective because the output of any input is subject to uncertainty and therefore one has to account for the uncertainty in some way in order to assess the overall robust performance.
\begin{remark}
	[RTS-R2 utility] Notice that one can rewrite an RTS-R2 utility in terms of a standard R2 utility: $R^{\textnormal{RTS}}_{f, \boldsymbol{\rho}}[X]=U[\cup_{\mathbf{x} \in X} \boldsymbol{\rho}[f(\mathbf{x}, \cdot)]]$. Therefore, to compute the RTS performance of a set of inputs, one just has to compute the corresponding R2 utility of the union risk-adjusted output set. In general, a similar result does not hold for the STR-R2 utility unless the operations commute (\cref{sec:commute}).
	\label{rem:rts_utility}
\end{remark}
\subsection{Robust inverted generational distance}
\label{sec:robust_igd}
The inverted generational distance (IGD) is a popular multi-objective performance metric that was initially introduced by \cite{coello2004m2aai}. Following \cite{schutze2012itec}, the IGD indicator $I^{\text{IGD}_{p,q}}: 2^{\mathbb{R}^M} \times B[\mathbb{R}^M] \rightarrow \mathbb{R}$ computes some notion of the average distance between a set of vectors $Y \subset \mathbb{R}^M$ and a finite set of ideal (or utopia) vectors $\Upsilon \in \mathbb{B}[\mathbb{R}^M]$,
\begin{equation}
	I^{\textnormal{IGD}_{p, q}}[Y, \Upsilon] = \biggl(\frac{1}{|\Upsilon|} \sum_{\boldsymbol{\upsilon} \in \Upsilon} \Bigl(\inf_{\mathbf{y} \in Y} ||\boldsymbol{\upsilon} - \mathbf{y}||_{L^p}\Bigr)^q\biggr)^{1/q},
	\label{eqn:igd}
\end{equation}
for some norms $p, q \geq 1$, where $\mathbb{B}[\mathbb{R}^M]:=\{Y \subseteq \mathbb{R}^M:|Y|<\infty\}$ denotes the space of finite sets of vectors. By construction, a smaller IGD indicator value is more desirable. Following \cite{tu2025sr}, we can rewrite the IGD indicator as a transformation of an R2 utility, namely the IGD utility
\begin{equation}
	U^{\textnormal{IGD}_{p, q}}[Y] := - (I^{\textnormal{IGD}_{p, q}}[Y, \Upsilon])^q = \frac{1}{|\Upsilon|} \sum_{\boldsymbol{\upsilon} \in \Upsilon} \sup_{\mathbf{y} \in Y} s^{\textnormal{IGD}_{p, q}}_{\boldsymbol{\upsilon}}(\mathbf{y}),
	\label{eqn:igd_utility}
\end{equation}
where $s^{\textnormal{IGD}_{p, q}}_{\boldsymbol{\upsilon}}(\mathbf{y}) = - ||\boldsymbol{\upsilon} - \mathbf{y}||_{L^p}^q$ is the IGD scalarisation function. Note that the total ordering imposed by the IGD utility \eqref{eqn:igd_utility} remains the same as the IGD indicator \eqref{eqn:igd} but with the directions reversed. That is, a smaller indicator value was preferred before, whereas now, a larger utility value is more desirable. By substituting this R2 utility set-up into the corresponding expressions for the robust R2 utilities \eqref{eqn:rts_utility} and \eqref{eqn:str_utility}, we obtain the following robust generalisations of the IGD utility:
\begin{align*}
	R_{f, \boldsymbol{\rho}}^{\textnormal{RTS-IGD}_{p, q}}[X] 
	&:= \frac{1}{|\Upsilon|} \sum_{\boldsymbol{\upsilon} \in \Upsilon} \sup_{\mathbf{x} \in X} \max_{\mathbf{y} \in \boldsymbol{\rho}[f(\mathbf{x}, \cdot)]} s^{\textnormal{IGD}_{p, q}}_{\boldsymbol{\upsilon}}(\mathbf{y}),
	\\
	R_{f, \rho}^{\textnormal{STR-IGD}_{p, q}}[X] 
	&:= \frac{1}{|\Upsilon|} \sum_{\boldsymbol{\upsilon} \in \Upsilon} \sup_{\mathbf{x} \in X} \rho[s^{\textnormal{IGD}_{p, q}}_{\boldsymbol{\upsilon}}(f(\mathbf{x}, \cdot))],
	\label{eqn:robust_igd_utility}
\end{align*}
which holds for any set of inputs $X \subseteq \mathbb{X}$ and norms $p, q \geq 1$. In other words, the RTS-IGD utility computes the average negative distance between the risk-adjusted output set and the set of utopia points, whilst the STR-IGD utility computes the average risk-adjusted negative distance between the distribution of possible objective vectors and the set of utopia points. On the left of \cref{fig:robust_performance_metrics}, we give a pictorial illustration of these different computations on a simple two-dimensional example.

\begin{remark}
	[Robust IGD+ utility] We can also generalise the above discussion for the corresponding IGD+ indicator described by \cite{ishibuchi2015emo}. This is a straightforward extension, where the IGD scalarisation function is replaced with its truncated version \citep{tu2025sr}.
\end{remark}
\subsection{Robust hypervolume indicator}
\label{sec:robust_hypervolume}
\begin{figure}
	\includegraphics[width=1\linewidth]{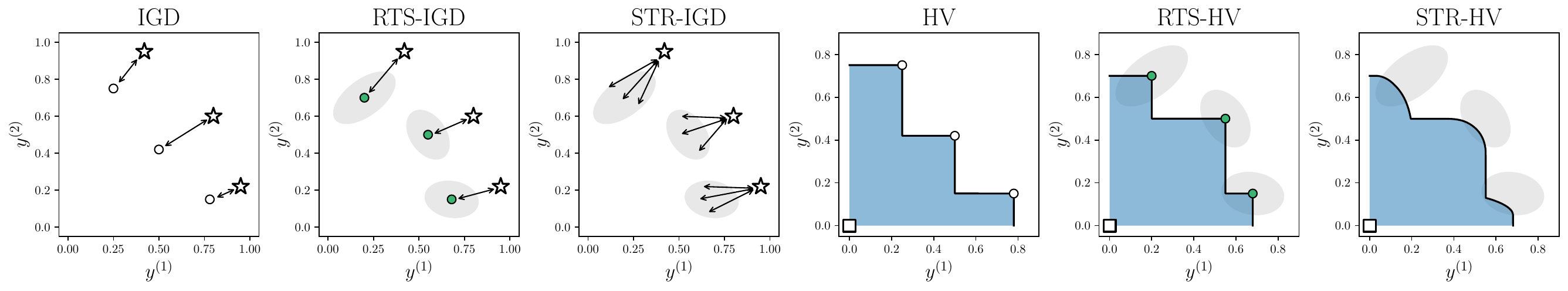}
	\centering
	\caption{An illustration of the standard and robust IGD (\cref{sec:robust_igd}) and hypervolume (\cref{sec:robust_hypervolume}) indicators on a two-dimensional example. In the plots, we illustrate the output associated with three inputs which are assumed to be either fixed for the standard (no uncertainty) examples, or uniformly distributed on some ellipsoid for the robust examples. For the robust utilities, we used the expectation functionals; for the IGD indicator, we use three utopia points (stars); and for the HV indicator, we set the reference point to the origin (square).}
	\label{fig:robust_performance_metrics}
\end{figure}

The hypervolume indicator (HV) is another popular multi-objective performance metric which was proposed by \cite{zitzler1998ppsn}. Formally, the hypervolume indicator $I^{\text{HV}}: 2^{\mathbb{R}^M} \times \mathbb{R}^M \rightarrow \mathbb{R}$ computes the volume that is enclosed between a reference vector $\boldsymbol{\eta} \in \mathbb{R}^M$ and the Pareto front surface of a set of vectors $Y \subset \mathbb{R}^M$, that is
\begin{equation}
	I^{\textnormal{HV}}[Y, \boldsymbol{\eta}] :=
	\int_{\mathbb{R}^M} \mathbbm{1}[\mathbf{z} \in \mathbb{D}_{\preceq}[\mathcal{Y}_{\boldsymbol{\eta}}^{\text{int}}[Y]] \cap \mathbb{D}_{\succsucc}[\{\boldsymbol{\eta}\}] ] d\mathbf{z}.
	\label{eqn:hypervolume}
\end{equation}
As identified in existing work \citep{shang2018pgecc,deng2019itec,zhang2020icml}, we can rewrite the hypervolume indicator as an R2 utility,
\begin{equation*}
	U_{\boldsymbol{\eta}}^{\textnormal{HV}}[Y] 
	:= I_{\textnormal{HV}}[Y, \boldsymbol{\eta}]
	= \mathbb{E}_{\boldsymbol{\lambda} \sim \textnormal{Uniform}(\mathcal{S}_+^{M-1})}
	\biggl[\sup_{\mathbf{y} \in Y} \tau^{\textnormal{HV}}(s^{\textnormal{Len}}_{(\boldsymbol{\eta}, \boldsymbol{\lambda})}(\mathbf{y})) \biggr],
\end{equation*}
where $\tau^{\text{HV}}(x) = c_M x^M$ is the hypervolume transformation function which depends on a constant $c_M = \pi^{M/2}/(2^M \Gamma(M/2 + 1))$ with $\Gamma(z) = \int_0^\infty t^{z-1} e^{-t} dt$ denoting the Gamma function. As we showcase below in \cref{prop:rts_hypervolume} and \cref{prop:str_hypervolume}, the hypervolume enclosed between a reference vector and the RTS and STR fronts can also naturally be written as an RTS-R2 and STR-R2 utility, respectively. We prove both of these results in \cref{app:proofs:prop:robust_hypervolume}. For some additional intuition, we also present a two-dimensional example of these different hypervolumes on the right of \cref{fig:robust_performance_metrics}.
\begin{proposition}
	[RTS hypervolume] Consider a bounded objective function $f: \mathbb{X} \times \Xi \rightarrow \mathbb{R}^M$, a bounded multivariate risk functional $\boldsymbol{\rho}$ and a reference vector $\boldsymbol{\eta} \in \mathbb{R}^M$. The hypervolume of the corresponding RTS front \eqref{eqn:rts_front} can be written as an RTS-R2 utility
	\begin{align*}
		R_{\boldsymbol{\eta}, f, \boldsymbol{\rho}}^{\textnormal{RTS-HV}}[X] 
		&:= \int_{\mathbb{R}^M} \mathbbm{1}[\mathbf{z} \in \mathbb{D}_{\preceq}[\mathcal{F}_{\boldsymbol{\eta}, f, \boldsymbol{\rho}}^{\textnormal{RTS}}[X]] \cap \mathbb{D}_{\succsucc}[\{\boldsymbol{\eta}\}]] d\mathbf{z}
		\\
		&= \mathbb{E}_{\boldsymbol{\lambda} \sim \textnormal{Uniform}(\mathcal{S}_+^{M-1})}
		\biggl[\sup_{\mathbf{x} \in X} \max_{\mathbf{y} \in \boldsymbol{\rho}[f(\mathbf{x}, \cdot)]} \tau^{\textnormal{HV}}(s^{\textnormal{Len}}_{(\boldsymbol{\eta}, \boldsymbol{\lambda})}(\mathbf{y})) \biggr].
	\end{align*}
	\label{prop:rts_hypervolume}
\end{proposition}
\begin{proposition}
	[STR hypervolume] Consider a bounded objective function $f: \mathbb{X} \times \Xi \rightarrow \mathbb{R}^M$, a partially coherent univariate risk functional $\rho$ and a reference vector $\boldsymbol{\eta} \in \mathbb{R}^M$. The hypervolume of the corresponding STR front \eqref{eqn:str_front} can be written as an STR-R2 utility
	\begin{align*}
		R_{\boldsymbol{\eta}, f, \rho}^{\textnormal{STR-HV}}[X] 
		&:= \int_{\mathbb{R}^M} \mathbbm{1}[\mathbf{z} \in \mathbb{D}_{\preceq}[\mathcal{F}_{\boldsymbol{\eta}, f, \rho}^{\textnormal{STR}}[X]] \cap \mathbb{D}_{\succsucc}[\{\boldsymbol{\eta}\}] ] d\mathbf{z}
		\\
		&= \mathbb{E}_{\boldsymbol{\lambda} \sim \textnormal{Uniform}(\mathcal{S}_+^{M-1})}
		\biggl[\sup_{\mathbf{x} \in X} \rho^{\textnormal{HV}}[s^{\textnormal{Len}}_{(\boldsymbol{\eta}, \boldsymbol{\lambda})}(f(\mathbf{x}, \cdot))] \biggr]
	\end{align*}
	where $\rho^{\textnormal{HV}} := \tau^{\textnormal{HV}} \circ \rho$.
	\label{prop:str_hypervolume}
\end{proposition}

\begin{remark}
	[Monte Carlo estimate] Computationally, the hypervolume indicators defined above can be estimated efficiently using a Monte Carlo estimate with $J>0$ samples of the positive unit vectors $\boldsymbol{\lambda}_j \sim \textnormal{Uniform}(\mathcal{S}_+^{M-1})$. For instance, one can estimate the STR hypervolume with the equation
	\begin{equation*}
		\hat{R}_{\boldsymbol{\eta}, f, \rho, J}^{\textnormal{STR-HV}}[X] 
		= \frac{1}{J} \sum_{j=1}^J \sup_{\mathbf{x} \in X} \rho^{\textnormal{HV}}[s^{\textnormal{Len}}_{(\boldsymbol{\eta}, \boldsymbol{\lambda}_j)}(f(\mathbf{x}, \cdot))],
	\end{equation*}
	for any set of inputs $X \subseteq \mathbb{X}$, bounded objective function $f: \mathbb{X} \times \Xi \rightarrow \mathbb{R}^M$, partially coherent univariate risk functional $\rho$ and reference vector $\boldsymbol{\eta} \in \mathbb{R}^M$. By a simple calculation, we see that the variance of this estimate can be uniformly upper bounded by a function dependent constant divided by the number of Monte Carlo samples, that is
	\begin{align*}
		\mathbb{V}\textnormal{ar}[\hat{R}_{\boldsymbol{\eta}, f, \rho, J}^{\textnormal{STR-HV}}[X]]
		&\leq \frac{1}{J} \mathbb{E}_{\boldsymbol{\lambda} \sim \textnormal{Uniform}(\mathcal{S}_+^{M-1})}
		\biggl[
		\sup_{\mathbf{x} \in X} \rho^{\textnormal{HV}}[s^{\textnormal{Len}}_{(\boldsymbol{\eta}, \boldsymbol{\lambda})}(f(\mathbf{x}, \cdot))]^2
		\biggr]
		\\
		&\leq \frac{c^2_M}{J} \sup_{(\mathbf{x}, \boldsymbol{\xi}) \in \mathbb{X} \times \Xi} ||f(\mathbf{x}, \boldsymbol{\xi}) - \boldsymbol{\eta}||^{2M}_{L^2} 
		\cdot \mathbbm{1}[f(\mathbf{x}, \boldsymbol{\xi}) \in \mathbb{D}_{\succsucc}[\{\boldsymbol{\eta}\}]].
	\end{align*}
	Similarly, one can also devise and bound the variance of the Monte Carlo estimate for the standard and RTS hypervolume indicators in the same way.
\end{remark}
\section{Numerical examples}
\label{sec:experiments}
We now present two numerical case studies in order to highlight the efficacy of the novel ideas that we have presented in this work. Firstly, in \cref{sec:cake}, we present the cake baking case study which illustrates a practical example where one has to decide between an RTS and STR solution. Secondly, in \cref{sec:rocket}, we present the robust rocket injector case study which showcases how one can use the concepts of a robust Pareto front and a robust performance metric in order to devise and assess the performance of some robust multi-objective optimisation algorithms. 
\subsection{Cake problem}
\label{sec:cake}
\begin{figure}
	\includegraphics[width=1\linewidth]{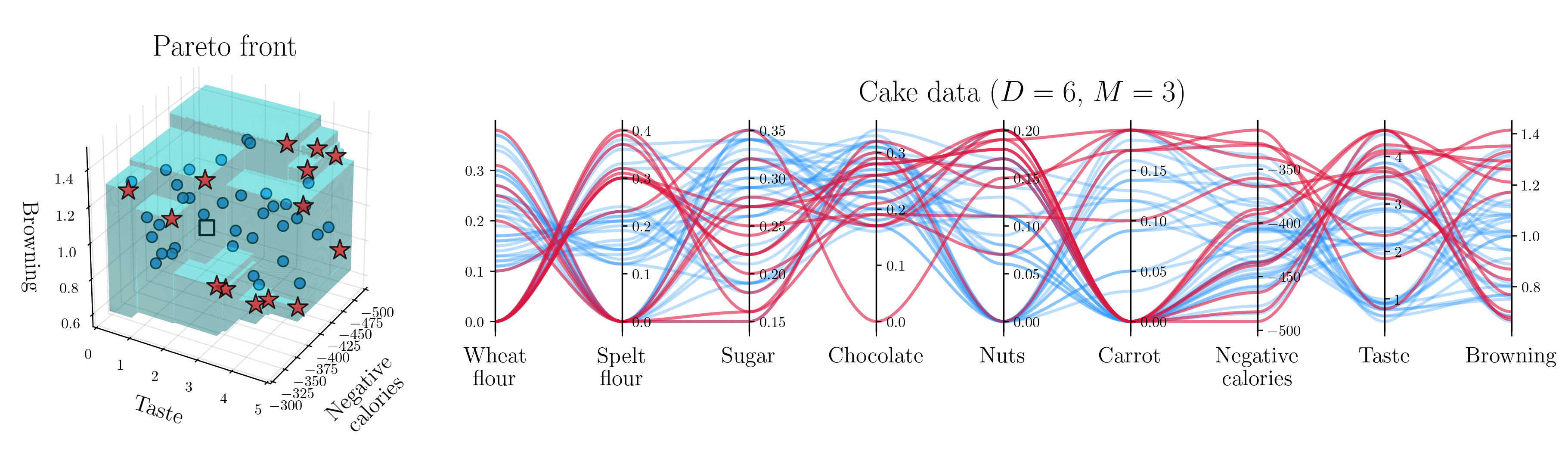}
	\centering
	\caption{An illustration of the cake data set. On the left, we plot the Pareto front surface of the $T=50$ observations. On the right, we draw the corresponding parallel coordinates plot over both the input and output space. The sub-optimal points are coloured in blue, whilst the Pareto optimal points are highlighted in red.}
	\label{fig:cake_parallel_coordinates}
\end{figure}

In this numerical example, we illustrate the utility of the RTS and STR approach through a real-world case study. This analysis is based on the cake data set $\mathcal{D}_T$, which is a public\footnote{The cake data set is available in the following Github repository: \href{https://github.com/basf/mopti}{\texttt{https://github.com/basf/mopti}}.} data set describing the results of $T=50$ cake baking experiments. Formally, this data set investigates how $M=3$ properties of a cake, described by the multi-objective function $g: \mathbb{X} \rightarrow \mathbb{R}^M$,
\begin{equation*}
	g(\mathbf{x}) = \small(-\textsc{Calories}(\mathbf{x}), \textsc{Taste}(\mathbf{x}), \textsc{Browning}(\mathbf{x})) \in \mathbb{R}^M,
\end{equation*}
changes when the relative composition of its $D=6$ main ingredients are varied. The inputs for this problem are the relative amount of ingredients $\mathbf{x} \in \mathbb{X}$, which lie in the $(D-1)$-dimensional simplex and satisfy the following inequality constraints:
\begin{align*}
	\mathbb{X} &= C_1 \cap C_2 \cap C_3 \cap C_4,
	\\
	C_1 &= \{\mathbf{x} \in \mathbb{R}_{\geq 0}^D: ||\mathbf{x}||_{L^1} = 1\},
	\\
	C_2 &= \{\mathbf{x} \in \mathbb{R}^D: x^{(1)} + x^{(2)} \in [0.2, 0.4] \},
	\\
	C_3 &= \{\mathbf{x} \in \mathbb{R}^D: x^{(3)} \in [0.15, 0.35] \},
	\\
	C_4 &= \{\mathbf{x} \in \mathbb{R}^D: x^{(4)} + x^{(5)} + x^{(6)} \geq 0.2 \}.
\end{align*}
Intuitively, these constraints are enforced in order to ensure that the resulting bakes lead to valid cakes. In words, condition one ($C_1$) is the simplex constraint; condition two ($C_2$) is a constraint on the amount of flour that is used within the bake; condition three ($C_3$) is a constraint on the amount of sugar that is used; and lastly condition four ($C_4$) is a constraint on the amount of chocolate, nuts and carrots that is used within the mix. In \cref{fig:cake_parallel_coordinates}, we illustrate this data set and highlight the corresponding $P=14$ observed Pareto optimal points.

\begin{figure}
	\includegraphics[width=1\linewidth]{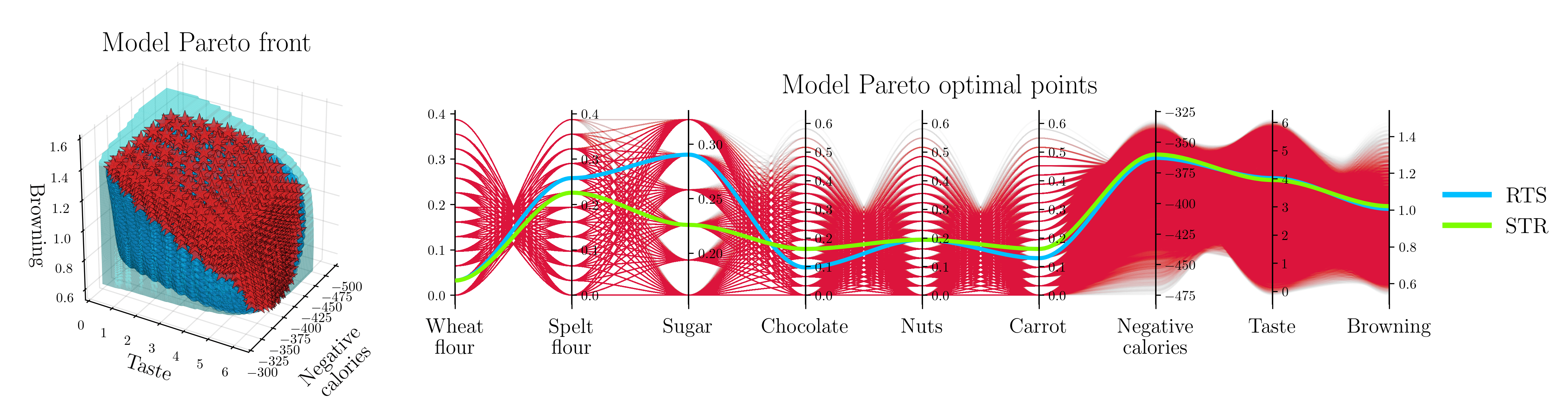}
	\centering
	\caption{An illustration of the mean function $\boldsymbol{\mu}_T: \mathbb{X} \rightarrow \mathbb{R}^M$ associated with a Gaussian process model that was trained on the cake data set. On the left, we plot the mean vectors over a discretisation of the input space and its corresponding Pareto front surface. On the right, we draw a parallel coordinates plot of the Pareto optimal vectors in red. In addition, we also highlight the RTS and STR solution for the target optimisation problem (\cref{sec:cake}) in light blue and green, respectively.}
	\label{fig:cake_model_parallel_coordinates}
\end{figure}

In practice, given a data set such as this, one wishes to solve the inverse problem. That is, to find a cake whose properties achieves a certain set of desirable values. As this data set only contains a finite number of experiments, we follow standard procedure and fit a probabilistic model on it. Intuitively, this model is used as a surrogate in order to get a prediction on the properties of a cake at any feasible set of ingredients. For convenience, in this example, we fit a Gaussian process model, which is a popular model that is often used in industrial applications---see for instance \cite{rasmussen2006} for a review. More precisely, we let $f: \mathbb{X} \times \Xi \rightarrow \mathbb{R}^M$, denote the corresponding posterior Gaussian process model that was fitted on the data points. That is, every instance of the uncertain parameter $\boldsymbol{\xi} \in \Xi$ corresponds to a different sample path of the Gaussian process model: 
\begin{equation*}
	(f(\mathbf{x}_1, \boldsymbol{\xi}), \dots, f(\mathbf{x}_N, \boldsymbol{\xi})) 
	\sim \mathcal{N}(\boldsymbol{\mu}_T(X), \boldsymbol{\Sigma}_T(X, X))
\end{equation*}
for any finite set of inputs $X = (\mathbf{x}_1,\dots, \mathbf{x}_N) \in \mathbb{X}^N$, where $\boldsymbol{\mu}_T: \mathbb{X} \rightarrow \mathbb{R}^M$ and $\boldsymbol{\Sigma}_T : \mathbb{X} \times \mathbb{X} \rightarrow \mathbb{R}^{M \times M}$ denotes the corresponding mean and covariance function of the posterior Gaussian process model. In \cref{fig:cake_model_parallel_coordinates}, we plot the posterior mean of this model over a discretisation of the input space with $|\hat{\mathbb{X}}| = 44835$. We also highlight the corresponding Pareto optimal vectors and Pareto front surface.

Suppose for instance that the decision maker is interested in identifying the ingredients that are required to make a cake whose properties achieve a target value of $\boldsymbol{\upsilon} = (-360, 4, 1) \in \mathbb{R}^M$. As motivated earlier in \cref{eg:distance_to_ideal}, one can recast this target optimisation problem as a robust scalarised optimisation problem. Specifically, one can use the expectation functionals (\cref{eg:expectation,eg:multivariate_expectation}) along with a distance-based scalarisation function. As a concrete example, we consider using the following weighted variant of the IGD scalarisation function:
\begin{equation*}
	s^{\text{WIGD}_{2,2}}_{(\boldsymbol{\upsilon}, \mathbf{w})}(\mathbf{y}) 
	:= - ||\mathbf{w}(\boldsymbol{\upsilon}-\mathbf{y})||^2_{L^2}
	:= - \sum_{m=1}^M \bigl(w^{(m)}(\upsilon^{(m)} - y^{(m)})\bigr)^2
\end{equation*}
for some weight vector $\mathbf{w} \in \Delta^{M-1}$. The introduction of this weight vector is mainly to ensure that each objective is normalised so that they have similar ranges. In particular, it ensures that the distance computation is not skewed towards the objectives with larger ranges such as the negative number of calories. For our example, we set the weight to be the normalised inverse of the estimated objective ranges: $\mathbf{w} = (0.004, 0.155, 0.841) \in \Delta^{M-1}$. Equipped with this set-up, the corresponding RTS and STR solutions for this target optimisation problem reduces down to the following two expressions:
\begin{align*}
	\mathbf{x}^*_{\text{RTS}} 
	&\in \argmin_{\mathbf{x} \in \mathbb{X}} ||\mathbf{w}(\boldsymbol{\upsilon}-\boldsymbol{\mu}_T(\mathbf{x}))||^2_{L^2},
	\\
	\mathbf{x}^*_{\text{STR}}
	&\in \argmin_{\mathbf{x} \in \mathbb{X}} \Bigl(||\mathbf{w}(\boldsymbol{\upsilon}-\boldsymbol{\mu}_T(\mathbf{x}))||^2_{L^2} + 
	||\mathbf{w} \sqrt{\textsc{Diagonal}(\boldsymbol{\Sigma}_T(\mathbf{x}, \mathbf{x}))}||^2_{L^2}
	\Bigr),
\end{align*}
where all operations on vectors are defined component-wise and the $\textsc{Diagonal}$ operation is used in order to extract the $M$-dimensional vector of variances from the square covariance matrix $\boldsymbol{\Sigma}_T(\mathbf{x}, \mathbf{x}) \in \mathbb{R}^{M\times M}$ for any input $\mathbf{x} \in \mathbb{X}$. Evidently, from these equations, we see that the difference between an RTS and STR approach is very similar to the classical bias-variance trade-off. The RTS problem is solely focussed on minimising the bias. That is, it tries to seek out the set of ingredients whose expected output is as close to the target as possible in terms of the weighted $L^2$-norm. In contrast, the STR approach also considers the effect of the variance in the scalarised problem as well. In particular, the STR problem contains an additional variance penalty term which penalises points whose weighted variances are large. In \cref{fig:cake_model_parallel_coordinates}, we plot the corresponding RTS and STR solution for the discretised problem:
\begin{align*}
	\mathbf{x}^*_{\text{RTS}} &= (0.0323, 0.2581, 0.2903, 0.0968, 0.1935, 0.1290),
	\\
	\mathbf{x}^*_{\text{STR}} &= (0.0323, 0.2258, 0.2258, 0.1613, 0.1935, 0.1613).
\end{align*}
Notably, the main difference between these two recipes lies in the relative amount of sugar that is used. The RTS approach suggests that we should have a recipe that uses more sugar if we are interested in getting an expected performance that is close to the target. In contrast, the STR approach suggests that we should substitute the additional sugar in favour of more chocolate and carrots. Conceptually, we should use the RTS solution over the STR solution in the aggregation scenario, where we assess the performance based on a collective aggregate. For instance, if we were making and selling boxes of little cupcakes, then this RTS solution would be more ideal because we would be interested in optimising the average quality in each box. In contrast, the STR solution becomes more useful in the scenario where the quality of each cake is assessed individually. For instance, if we were selling many individual cakes, then we should favour the STR solution as the quality of each cake should be assessed independently in this setting.
\subsection{Robust rocket injector problem}
\label{sec:rocket}
In this numerical example, we study a distributionally robust variant of the rocket injector problem. This problem was initially formulated in the work by \cite{vaidyanathan20034asmae} and was later implemented by \cite{tanabe2020asc}. Formally, this problem studies the optimisation of $M=3$ (scaled) objectives relating to the performance and longevity of a rocket injector design $f:\mathbb{X} \times \Xi \rightarrow \mathbb{R}^M$,
\begin{equation*}
	f(\mathbf{x}, \boldsymbol{\xi}) = \small(
	-\textsc{MaxFaceTemp}(\mathbf{x}, \boldsymbol{\xi}), 
	-\textsc{CombustionLength}(\mathbf{x}, \boldsymbol{\xi}),
	-\textsc{MaxOxTipTemp}(\mathbf{x}, \boldsymbol{\xi})
	) \in \mathbb{R}^M,
\end{equation*}
subject to changes in $D+W=4$ geometric features $(\mathbf{x}, \boldsymbol{\xi}) \in \mathbb{X} \times \Xi$.  Unlike the original work, we consider a robust variant of this problem in which only a subset of the geometric features can be adjusted by the decision maker. The other components are assumed to be determined by the supply manufacturer who provides a base construction which is subject to some uncertainty. In particular, we suppose that we can control $D=2$ geometric features lying in the normalised space $\mathbb{X} = [0, 1]^D$, whilst the manufacturer can control the other $W=2$ pertinent features lying in the normalised space $\Xi = [0, 1]^W$. Explicitly, the controllable features are the ones affecting the Hydrogen supply, whilst the uncontrollable features affect the Oxygen supply:
\begin{align*}
	x^{(1)} &: \text{Normalised angle at which the $H_2$ is directed toward the oxidiser},
	\\
	x^{(2)} &: \text{Normalised change in the $H_2$ flow area from the baseline},
	\\
	\xi^{(1)} &: \text{Normalised change in the $O_2$ flow area from the baseline},
	\\
	\xi^{(2)} &: \text{Normalised oxidiser post-tip thickness}.
\end{align*}
For the manufacturer's distribution, we assume a uniform distribution $p(\boldsymbol{\xi}) \sim \text{Uniform}(\Xi)$ over the uncertain parameters $\boldsymbol{\xi} \in \Xi$. For our robust objective, we consider an STR approach where the goal is to identify the STR robust designs whose output distribution lies on the CVaR Pareto front surface at level $\alpha=0.9$. Equivalently, we formulate this problem as a robust utility optimisation problem where we are interested in optimising the STR hypervolume (\cref{prop:str_hypervolume}) of the identified CVaR front
\begin{equation*}
	\max_{X \subseteq \mathbb{X}, |X| \leq T} R_{\boldsymbol{\eta}, f, \rho^{\text{CVaR}_{\alpha}}}^{\textnormal{STR-HV}}[X],
\end{equation*}
where the reference vector $\boldsymbol{\eta} = (-1.05, -1.30, -1.20) \in \mathbb{R}^M$ is set to a lower bound of the objective function and $T > 0$ is a cardinality constraint. Note that the inclusion of this cardinality constraint is crucial because it encodes the fact that we are interested in finding a finite set of high quality points as opposed to a large set of potentially variable quality points. In fact, without this constraint, one can trivially solve this problem by proposing the whole input space $X = \mathbb{X}$ because the hypervolume is monotonic and therefore there is no loss in utility with adding more points.

\begin{figure}
	\includegraphics[width=1\linewidth]{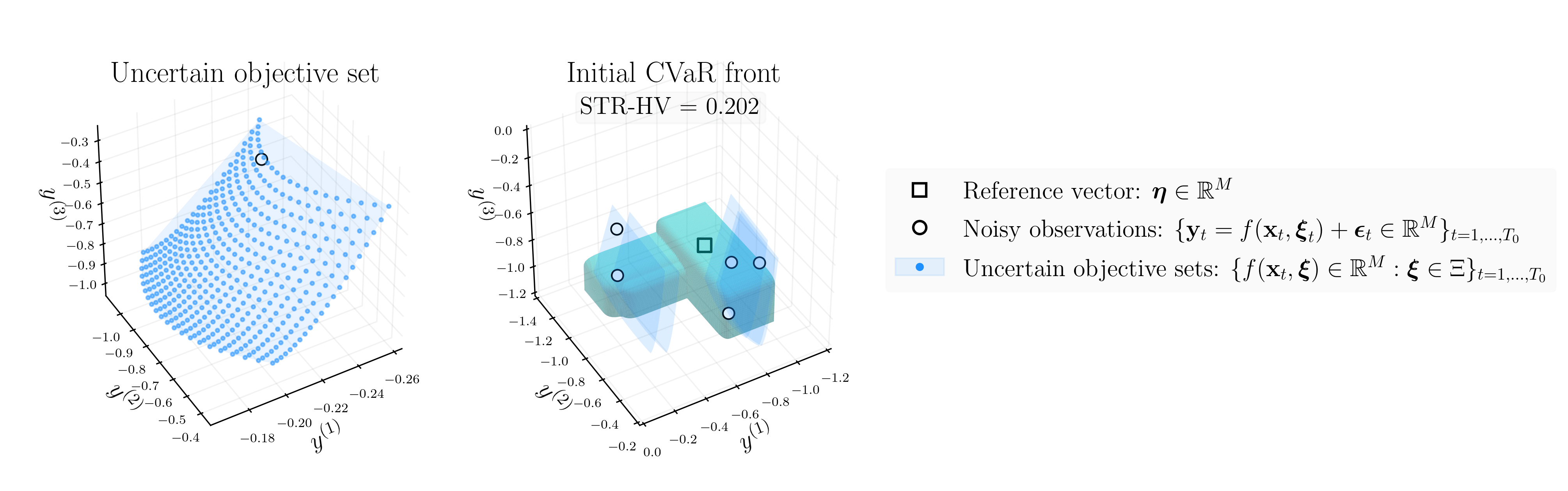}
	\centering
	\caption{This plot illustrates the initial CVaR front at level $\alpha=0.9$ for the $T_0=5$ initial designs under the STR approach. On the left, we plot the uncertainty objective set associated with a single input location. On the right, we plot the CVaR front associated with the initial set of inputs.}
	\label{fig:rocket_cvar_front_initial}
\end{figure}

To emulate a real-world scenario, we treat the objective function $f$ as a black-box function that can only be evaluated sparingly at a few carefully selected design locations $T=65$. More precisely, at any time $t > 0$, the decision maker is given the opportunity to see a noisy function evaluation at any chosen design $\mathbf{x}_t \in \mathbb{X}$, that is $\mathbf{y}_t = f(\mathbf{x}_t, \boldsymbol{\xi}_t) + \boldsymbol{\epsilon}_t \in \mathbb{R}^M$, where $\boldsymbol{\xi}_t \sim p(\boldsymbol{\xi})$ is the uncertain parameter that was independently sampled by the manufacturer and $\boldsymbol{\epsilon}_t \sim \mathcal{N}(\mathbf{0}_M, \textsc{Diag}(\boldsymbol{\sigma}^2))$ is the observation noise. In our example, we consider a zero-mean Gaussian observation noise whose standard deviation $\boldsymbol{\sigma} = (0.01, 0.01, 0.01)$ is set to $1\%$ of the estimated scaled objective range. To initialise this problem, we randomly sampled $T_0=5$ initial points from the controllable input space. In \cref{fig:rocket_cvar_front_initial}, we illustrate these initial points in the objective space and plot the corresponding CVaR front under the STR approach. Notably, to compute this CVaR front, we require knowledge of the true objective function at each uncertain parameter, which is typically not available in a real-world set-up.

To approximately solve this robust multi-objective optimisation problem, we require an acquisition procedure in order to sequentially select the inputs to evaluate at. Note that we cannot deploy the STR approach directly for this problem because we do not have access to the true objective function. Instead, we propose using a straightforward extension of the standard hypervolume improvement acquisition function \citep{emmerich2006itec}. More specifically, at each time $t>0$ we propose choosing the input $\mathbf{x}_t \in \mathbb{X}$ which maximises the STR hypervolume improvement under the latest surrogate model
\begin{equation*}
	\mathbf{x}_t \in \argmax_{\mathbf{x} \in \mathbb{X}} R_{\hat{\boldsymbol{\eta}}_{t-1}, \hat{f}_{t-1}, \rho^{\text{CVaR}_{\alpha}}}^{\textnormal{STR-HV}}[X_{t-1} \cup \{\mathbf{x}\}] - R_{\hat{\boldsymbol{\eta}}_t, \hat{f}_{t-1}, \rho^{\text{CVaR}_{\alpha}}}^{\textnormal{STR-HV}}[X_{t-1}],
\end{equation*}
where $X_{t-1} = \{\mathbf{x}_n\}_{n\leq t-1}$ denotes the set of previously queried inputs, whilst $\hat{\boldsymbol{\eta}}_t \in \mathbb{R}^M$ and $\hat{f}_{t-1}: \mathbb{X} \times \Xi \rightarrow \mathbb{R}^M$ denotes the latest estimate of the reference vector and objective function, respectively. In our example, we set the latest reference vector to be an over-estimation of the lower bound $\eta_{t-1}^{(m)} = \min_{n < t} y_n^{(m)} - \kappa (\max_{n < t} y_n^{(m)} - \min_{n < t} y_n^{(m)})$ with $\kappa = 0.1$ for $m=1,\dots,M$. For the surrogate model, we used an upper confidence bound of a Gaussian process model which was trained on the previous observations, that is $\hat{f}^{\text{UCB}}_{t-1}(\mathbf{x}, \boldsymbol{\xi}) = \boldsymbol{\mu}_{t-1}(\mathbf{x}, \boldsymbol{\xi}) + \beta \boldsymbol{\sigma}_{t-1}(\mathbf{x}, \boldsymbol{\xi})$, where $\beta = 2.0$ is the trade-off parameter and the vectors $\boldsymbol{\mu}_{t-1}(\mathbf{x}, \boldsymbol{\xi}) \in \mathbb{R}^M$ and $\boldsymbol{\sigma}_{t-1}(\mathbf{x}, \boldsymbol{\xi}) \in \mathbb{R}_{\geq 0}^M$ are the mean and standard deviation of the model at the input pair $(\mathbf{x}, \boldsymbol{\xi}) \in \mathbb{X} \times \Xi$. We refer to this overall acquisition strategy as the STR-HVI-UCB approach. 
\begin{figure}
	\includegraphics[width=1\linewidth]{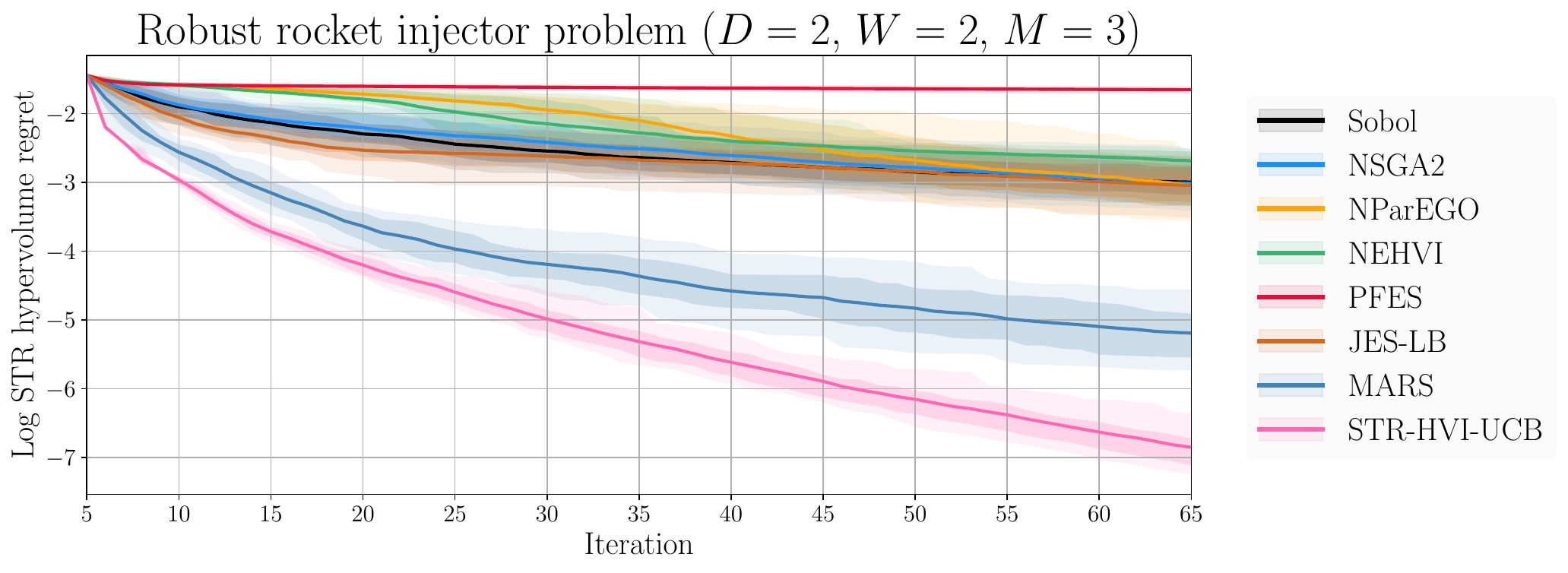}
	\centering
	\caption{A numerical comparison between the different multi-objective algorithms on the robust rocket injector problem over $100$ independent runs. For each method, we plot the mean log regret at every iteration along with an uncertainty estimate. Specifically, we highlight the band between the $0.25$ to $0.75$ and $0.1$ to $0.9$ quantiles with a dark and light shade, respectively.}
	\label{fig:rocket_cvar_front_results}
\end{figure}

\begin{figure}
	\includegraphics[width=1\linewidth]{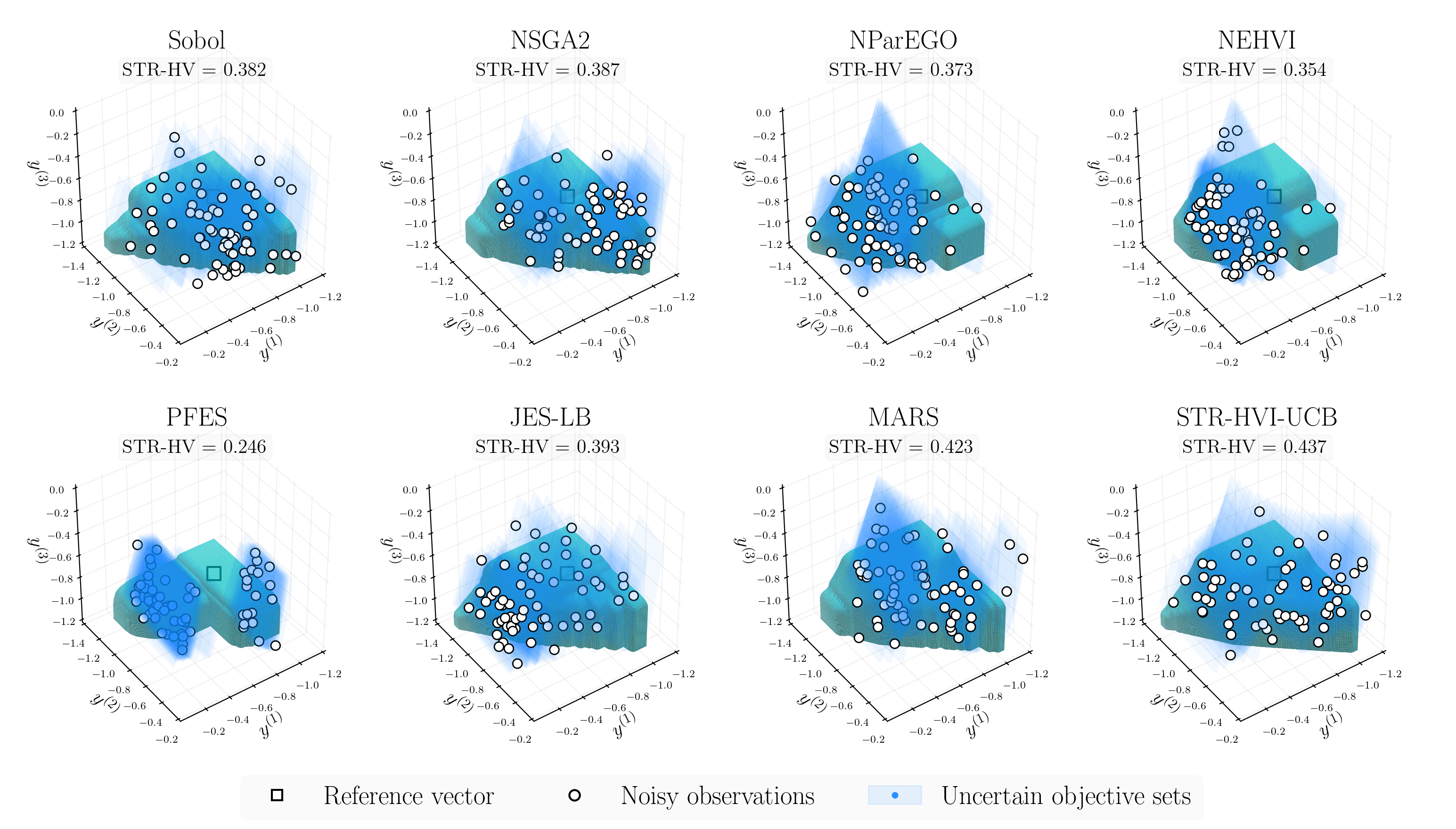}
	\centering
	\caption{An illustration of the final CVaR fronts identified by the different multi-objective optimisation algorithms on the robust rocket injector problem for one experimental run with $T=65$ function evaluations. The maximum hypervolume on this example is around $R_{\boldsymbol{\eta}, f, \rho^{\textnormal{CVaR}_{\alpha}}}^{\textnormal{STR-HV}}[\mathbb{X}] \approx 0.4382$.}
	\label{fig:rocket_cvar_fronts_example}
\end{figure}

In \cref{fig:rocket_cvar_front_results}, we present the results of our numerical comparison in which we compared our proposed method versus the existing baselines over $100$ independent runs on a discretised version of the robust rocket injector problem. Namely, we compare our method against non-robust\footnote{These approaches are non-robust in the sense that they can condition on the observed data set $\mathcal{D}_{t-1} = \{(\mathbf{x}_n, \boldsymbol{\xi}_n, \mathbf{y}_n)\}_{n\leq t-1}$, at every time $t$, but they do not consider the uncertainty present in the uncontrollable variable $\boldsymbol{\xi} \in \Xi$.} approaches such as: the random search (Sobol), NSGA2 \citep{deb2002itec}, NParEGO \citep{knowles2006itec,daulton2021anips}, NEHVI \citep{emmerich2006itec,daulton2021anips}, PFES \citep{suzuki2020icml} and the JES-LB \citep{tu2022anips} algorithms. We also compare our approach against the MARS algorithm \citep{daulton2022icml}, which is a robust algorithm that takes an STR approach using the Chebyshev scalarisation function and the univariate VaR risk functional. Unsurprisingly we see that the robust methods perform better than the non-robust methods. Moreover, our method performs the best overall under the chosen utility. For some visual intuition, we have included \cref{fig:rocket_cvar_fronts_example} where we plot the final CVaR front obtained by each algorithm under the same random seed. Visually we see that most of the variation between these algorithms transpire around the first and second objectives. The existing approaches regularly struggle to find the inputs which robustly maximises the CVaR front over these two objectives. In contrast, our robust approach is much more consistent and is often very quick in identifying the difficult inputs which robustly maximise the whole CVaR front.
\section{Conclusion}
\label{sec:conclusion}
In this work, we studied the robust multi-objective optimisation problem from a computational perspective. We identified that the majority of existing robust multi-objective optimisation algorithms, in one way or another, rely on two key operations: robustification and scalarisation. These two procedures are in general not commutative and therefore the order that they are applied in has an effect on the resulting solutions that are identified and the final decisions that are made. As part of our work, we analysed the effect of this ordering from both a philosophical and a methodological perspective. Notably, among other things, we proposed the novel concepts of a robust Pareto front surface and a robust performance metric. We then demonstrated the practical utility of these new ideas via some insightful numerical experiments.

In closing, the key ideas and methodology which we presented in this work are very general and encapsulates many existing ideas and algorithms into a neat computational framework. Future work could focus on refining these ideas further for specific problems and case studies. For instance, practitioners in industry are often interested in finding and using optimisation algorithms which are both numerically performant and theoretically justified. Many such algorithms and results has been developed for the robust single-objective optimisation problem over the past decades \citep{ben-tal2009,bertsimas2011sra,rahimian2019a}. Extending these computational algorithms and theoretical results to robust multi-objective setting is an interesting and practically useful direction for future research.

\section*{Acknowledgements}
Ben Tu was supported by the EPSRC StatML CDT programme EP/S023151/1 and BASF SE, Ludwigshafen am Rhein. Nikolas Kantas was partially funded by JPMorgan Chase \& Co. under J.P. Morgan A.I. Faculty Research Awards 2021.
\newpage
\appendix
\section{Proofs}
\label{app:proofs}

\subsection{Proof of \cref{thm:chebyshev_scalarisation}}
The proof of this result is presented by \citet[Part 2, Theorem 3.4.5]{miettinen1998}

\begin{flushright}
	$\blacksquare$
\end{flushright}
\subsection{Proof of \cref{lemma:bounded_univariate_risk}}
\label{app:proofs:lemma:bounded_univariate_risk}
Consider any partially coherent univariate risk functional $\rho$ and bounded function $h: \Xi \rightarrow \mathbb{R}$. We will show that $\rho[h]$ is bounded and therefore $\rho$ is a bounded univariate risk functional. As $h$ is bounded, there exists constants $l, u \in \mathbb{R}$ such that $l \leq h(\boldsymbol{\xi}) \leq u$ for all $\boldsymbol{\xi} \in \Xi$. Let $r \in \mathbb{R}$ be any constant and define $c_r: \Xi \rightarrow \mathbb{R}$ to be the corresponding constant function satisfying $c_r(\boldsymbol{\xi}) = r$ for all $\boldsymbol{\xi} \in \Xi$. Then by the normalised and translation equivariance property, we have that $\rho[c_r] = r$. By the monotonicity property, we have that $\rho[c_l] \leq \rho[h] \leq \rho[c_u]$. Together, these two results imply that $\rho[h]$ is bounded since $l \leq \rho[h] \leq u$.
\begin{flushright}
	$\blacksquare$
\end{flushright}
\subsection{Proof of \cref{prop:robust_chebyshev_scalarisation}}
\label{app:proofs:prop:robust_chebyshev_scalarisation}
Consider any weakly union robust input $\mathbf{z} \in \mathbb{X}$ for the union robust optimisation problem \eqref{eqn:robust_moo_union}. By the definition of weak union robustness, there must exist an output vector $\mathbf{y}^* \in \boldsymbol{\rho}[f(\mathbf{z}, \cdot)]$ which is weakly optimal. We will now proceed to show that the input $\mathbf{z}$ must be a maximiser to the RTS length scalarised optimisation problem \eqref{eqn:rts_length_problem} with the positive direction $\boldsymbol{\lambda} = (\mathbf{y}^* - \boldsymbol{\eta}) / ||\mathbf{y}^* - \boldsymbol{\eta}||_{L^2} \in \mathcal{S}^{M-1}_+$. Suppose for a contradiction that there exists a distinctive input $\mathbf{x} \in \mathbb{X}$ which achieves a greater RTS length scalarised value than the input $\mathbf{z}$ along this positive direction. This implies that
\begin{align*}
	\max_{\mathbf{y} \in \boldsymbol{\rho}[f(\mathbf{x}, \cdot)]} \min_{m=1,\dots,M} \frac{y^{(m)} - \eta^{(m)}}{\lambda^{(m)}}
	&> 
	\max_{\mathbf{y} \in \boldsymbol{\rho}[f(\mathbf{z}, \cdot)]} \min_{m=1,\dots,M} \frac{y^{(m)} - \eta^{(m)}}{\lambda^{(m)}}
	= ||\mathbf{y}^* - \boldsymbol{\eta}||_{L^2}
\end{align*}
due to the weak optimality of $\mathbf{y}^*$ and the strong monotonicity of the length scalarisation function \citep[Lemma A.2]{tu2024aa}. This further implies that
\begin{align*}
	\max_{\mathbf{y} \in \boldsymbol{\rho}[f(\mathbf{x}, \cdot)]} \min_{m=1,\dots,M} \frac{y^{(m)} - \eta^{(m)}}{y^{*(m)} - \eta^{(m)}} > 1,
\end{align*}
which means that there exists a vector $\mathbf{y} \in \boldsymbol{\rho}[f(\mathbf{x}, \cdot)]$ that strongly dominates $\mathbf{y}^*$. This contradicts the weak optimality of $\mathbf{y}^*$ and therefore concludes the proof.

\begin{flushright}
	$\blacksquare$
\end{flushright}
\subsection{Proof of \cref{thm:polar_parameterisation}}
\label{app:proofs:thm:polar_parameterisation}
The proof of this result is presented by \citet[Appendix A.1]{tu2024aa}.
\begin{flushright}
	$\blacksquare$
\end{flushright}

\subsection{Pareto front conditions}
\label{app:pareto_front_conditions}
In order to verify whether a polar surface is a Pareto front surface, we will rely on the key result by \citet[Proposition 3.1]{tu2024aa}, which we recall below in \cref{prop:pareto_front_conditions}.
\begin{proposition}
	[Pareto front conditions] \citep[Proposition 3.1]{tu2024aa} Consider a polar surface $A \in \mathbb{L}_{\boldsymbol{\eta}}$, then this set is a Pareto front surface, that is $A = \mathcal{Y}_{\boldsymbol{\eta}}^{\textnormal{int}}[A]$, if and only if the following two conditions holds:
	\begin{enumerate}[label=C\arabic*., ref=C\arabic*]
		\item The positive lengths condition: $\ell_{\boldsymbol{\eta}, \boldsymbol{\lambda}}[A] > 0$ for all $\boldsymbol{\lambda} \in \mathcal{S}_+^{M-1}$.
		\label{eqn:pareto_condition_1}
		\item The maximum ratio condition: $\max_{m=1,\dots,M} \frac{\ell_{\boldsymbol{\eta}, \boldsymbol{\lambda}}[A] \lambda^{(m)}}{\ell_{\boldsymbol{\eta}, \boldsymbol{\upsilon}}[A] \upsilon^{(m)}} \geq 1$ for all $\boldsymbol{\lambda}, \boldsymbol{\upsilon} \in \mathcal{S}_+^{M-1}$.
		\label{eqn:pareto_condition_2}
	\end{enumerate}
	\label{prop:pareto_front_conditions}
\end{proposition}
\subsection{Proof of \cref{prop:rho_front_statistics}}
\label{app:proofs:prop:rho_front_statistics}
Consider the Pareto front surface statistic defined in \eqref{eqn:rho_front_statistics}. Firstly, we will show that this set is indeed a polar surface. By definition, to establish that this set is a polar surface, we need to establish that the projected lengths
\begin{equation}
	\rho[\ell_{\boldsymbol{\eta}, \boldsymbol{\lambda}}[Y_{\boldsymbol{\eta}, f}^*(\cdot)]] 
	= \rho\Bigl[\sup_{\mathbf{x} \in \mathbb{X}} s^{\textnormal{Len}}_{(\boldsymbol{\eta}, \boldsymbol{\lambda})}(f(\mathbf{x}, \cdot)) \Bigr]
	\label{eqn:rho_length}
\end{equation}
are non-negative and bounded along every positive direction $\boldsymbol{\lambda} \in \mathcal{S}_+^{M-1}$. This result follows directly from the definition of the length scalarisation function \eqref{eqn:length_scalarisation} and the boundedness of the partially coherent univariate risk functional (\cref{lemma:bounded_univariate_risk}). We will now show that this polar surface is either a Pareto front surface or the degenerate singleton set. The former setting (Case 1) happens whenever there exists no positive direction where the projected length is zero. The latter setting (Case 2) happens whenever there exists at least one positive direction where the projected length is zero. Below, we consider these two cases separately.

\paragraph{Case 1:} Suppose that the projected length \eqref{eqn:rho_length} is positive along all positive directions. Then the projected length condition \eqref{eqn:pareto_condition_1} holds. We will now show that the maximum ratio condition \eqref{eqn:pareto_condition_2} also holds and therefore the polar surface \eqref{eqn:rho_front_statistics} is a Pareto front surface. Note that for any $\boldsymbol{\xi} \in \Xi$, the set $Y_{\boldsymbol{\eta}, f}^*(\boldsymbol{\xi})$ \eqref{eqn:stochastic_pareto_front} is either a valid Pareto front surface or a degenerate singleton set. In either case, it satisfies the maximum ratio condition
\begin{equation*}
	\ell_{\boldsymbol{\eta}, \boldsymbol{\lambda}}[Y_{\boldsymbol{\eta}, f}^*(\boldsymbol{\xi})] 
	\max_{m=1,\dots,M} \frac{\lambda^{(m)}}{\upsilon^{(m)}} 
	\geq \ell_{\boldsymbol{\eta}, \boldsymbol{\upsilon}}[Y_{\boldsymbol{\eta}, f}^*(\boldsymbol{\xi})]
\end{equation*}
for all $\boldsymbol{\lambda}, \boldsymbol{\upsilon} \in \mathcal{S}^{M-1}_+$ and $\boldsymbol{\xi} \in \Xi$. By the monotonicity and positive homogeneity of $\rho$, we get the desired result:
\begin{equation*}
	\rho[\ell_{\boldsymbol{\eta}, \boldsymbol{\lambda}}[Y_{\boldsymbol{\eta}, f}^*(\cdot)]]
	\max_{m=1,\dots,M} \frac{\lambda^{(m)}}{\upsilon^{(m)}} 
	\geq \rho[\ell_{\boldsymbol{\eta}, \boldsymbol{\upsilon}}[Y_{\boldsymbol{\eta}, f}^*(\cdot)]].
\end{equation*}
\paragraph{Case 2:} Now consider the setting where there exists one direction $\boldsymbol{\lambda} \in \mathcal{S}^{M-1}_+$ in which the projected length \eqref{eqn:rho_length} is zero. This implies that $\max_{\mathbf{x} \in \mathbb{X}} s^{\textnormal{Len}}_{(\boldsymbol{\eta}, \boldsymbol{\lambda})}(f(\mathbf{x}, \boldsymbol{\xi}))$ is zero on the support of $\rho$. As the univariate risk functional is partially coherent and the length scalarisation function is strongly monotonic \citep[Lemma A.2]{tu2024aa}, this implies that all other directions must also have a projected length of zero. Therefore, in this setting, the polar surface \eqref{eqn:rho_front_statistics} reduces down to the degenerate singleton set. 
\begin{flushright}
	$\blacksquare$
\end{flushright}
\subsection{Proof of \cref{prop:rts_front}}
\label{app:proofs:prop:rts_front}
Consider the RTS front defined in \eqref{eqn:rts_front}. Firstly, we will show that this set is indeed a polar surface. By definition, to establish that this set is a polar surface, we need to establish that the RTS lengths
\begin{equation}
	\sup_{\mathbf{x} \in X} \max_{\mathbf{y} \in \boldsymbol{\rho}[f(\mathbf{x}, \cdot)]} s^{\textnormal{Len}}_{(\boldsymbol{\eta}, \boldsymbol{\lambda})}(\mathbf{y})
	\label{eqn:rts_length}
\end{equation}
are non-negative and bounded along every positive direction $\boldsymbol{\lambda} \in \mathcal{S}_+^{M-1}$. This follows simply from the definition of the length scalarisation function \eqref{eqn:length_scalarisation} and the boundedness of both the objective function and multivariate risk functional. We will now show that the RTS front is either a Pareto front surface or the degenerate singleton set. These two cases are determined by whether the Pareto front surface $\mathcal{Y}_{\boldsymbol{\eta}}^{\textnormal{int}}[\cup_{\mathbf{x} \in X} \boldsymbol{\rho}[f(\mathbf{x}, \cdot)]]$ is non-empty (Case 1) or empty (Case 2), respectively. We will consider both of these separate cases below.

\paragraph{Case 1:} Suppose that the Pareto front surface of interest is non-empty. We can then invoke \cref{thm:polar_parameterisation} on the Pareto front surfaces $\mathcal{Y}_{\boldsymbol{\eta}}^{\textnormal{int}}[\cup_{\mathbf{x} \in X} \boldsymbol{\rho}[f(\mathbf{x}, \cdot)]]$ and $\mathcal{Y}_{\boldsymbol{\eta}}^{\textnormal{int}}[
\cup_{\mathbf{x} \in X} \mathcal{Y}_{\boldsymbol{\eta}}^{\textnormal{int}}[\boldsymbol{\rho}[f(\mathbf{x}, \cdot)]]]$. We see very quickly that both of these expressions reduce down to the RTS front \eqref{eqn:rts_front} and therefore the latter set is indeed a Pareto front surface in this setting.

\paragraph{Case 2:} Suppose that the Pareto front surface of interest is empty. This implies that all of the vectors in the set $\cup_{\mathbf{x} \in X} \boldsymbol{\rho}[f(\mathbf{x}, \cdot)]$ are weakly dominated by the reference vector. Consequently, this means that all of the RTS lengths \eqref{eqn:rts_length} are zero, which implies degeneracy: $\mathcal{F}_{\boldsymbol{\eta}, f, \boldsymbol{\rho}}^{\textnormal{RTS}}[X] = \{\boldsymbol{\eta}\} \in \mathbb{L}_{\boldsymbol{\eta}}$.
\begin{flushright}
	$\blacksquare$
\end{flushright}
\subsection{Proof of \cref{prop:str_front}}
\label{app:proofs:prop:str_front}
Consider the STR front defined in \eqref{eqn:str_front}. Firstly, we will show that this set is indeed a polar surface. By definition, to establish that this set is a polar surface, we need to establish that the STR lengths
\begin{equation}
	\sup_{\mathbf{x} \in X} \rho[s^{\textnormal{Len}}_{(\boldsymbol{\eta}, \boldsymbol{\lambda})}(f(\mathbf{x}, \cdot))]
	\label{eqn:str_length}
\end{equation}
are non-negative and bounded along every positive direction $\boldsymbol{\lambda} \in \mathcal{S}_+^{M-1}$. This follows simply from the definition of the length scalarisation function \eqref{eqn:length_scalarisation}, the boundedness of the objective function and the boundedness of the partially coherent univariate risk functional (\cref{lemma:bounded_univariate_risk}). We will now show that the STR front is either a Pareto front surface or the degenerate singleton set. These two cases are determined by whether the Pareto front surface $\mathcal{Y}_{\boldsymbol{\eta}}^{\textnormal{int}}[\cup_{\mathbf{x} \in X} \rho[\mathcal{Y}_{\boldsymbol{\eta}}^{\textnormal{int}}[f(\mathbf{x}, \cdot)]]]$ is non-empty (Case 1) or empty (Case 2), respectively. We will first recall a useful result before considering both of these separate cases below. Recall that the corresponding risk statistic of the individual Pareto front surface \eqref{eqn:rho_front_statistics} is defined by the set
\begin{align*}
	\begin{split}
		\rho[\mathcal{Y}_{\boldsymbol{\eta}}^{\textnormal{int}}[f(\mathbf{x}, \cdot)]]
		:= \{
		\boldsymbol{\eta} + 
		\rho[s^{\textnormal{Len}}_{(\boldsymbol{\eta}, \boldsymbol{\lambda})}(f(\mathbf{x}, \cdot))]
		\boldsymbol{\lambda}
		\in \mathbb{R}^M: \boldsymbol{\lambda} \in \mathcal{S}_+^{M-1} \} \in \mathbb{L}_{\boldsymbol{\eta}},
	\end{split}
\end{align*}
for every input $\mathbf{x} \in X$. As $\rho$ is a partially coherent univariate risk functional, the following equivalence holds along every positive direction $\boldsymbol{\lambda} \in \mathcal{S}_+^{M-1}$:
\begin{equation}
	s^{\textnormal{Len}}_{(\boldsymbol{\eta}, \boldsymbol{\lambda})}(\rho[\mathcal{Y}_{\boldsymbol{\eta}}^{\textnormal{int}}[f(\mathbf{x}, \cdot)]])
	= \rho[s^{\textnormal{Len}}_{(\boldsymbol{\eta}, \boldsymbol{\lambda})}(f(\mathbf{x}, \cdot))].
	\label{eqn:length_equivalence}
\end{equation}
This equivalence follows directly from the polar parameterised form of the risk statistic and the definition of the length scalarisation function. In essence it says that the point within the risk statistic which achieves the largest projected length along the line $\{\boldsymbol{\eta}+t\boldsymbol{\lambda} \in \mathbb{R}^M: t \in \mathbb{R}\}$ is indeed the vector $\boldsymbol{\eta} + \rho[s^{\textnormal{Len}}_{(\boldsymbol{\eta}, \boldsymbol{\lambda})}(f(\mathbf{x}, \cdot))] \boldsymbol{\lambda} \in \mathbb{R}^M$. 

\paragraph{Case 1:} Suppose that the Pareto front surface of interest is non-empty. We can then invoke \cref{thm:polar_parameterisation} on this Pareto front surface. By using the equivalence in \eqref{eqn:length_equivalence}, we see very quickly that this polar representation reduces down to the STR front \eqref{eqn:str_front} and therefore the latter set is indeed a valid Pareto front surface in this setting.

\paragraph{Case 2:} Suppose that the Pareto front surface of interest is empty. This implies that the sets $\rho[\mathcal{Y}_{\boldsymbol{\eta}}^{\textnormal{int}}[f(\mathbf{x}, \cdot)]]$ are degenerate for all inputs $\mathbf{x} \in X$. Consequently, by the equivalence in \eqref{eqn:length_equivalence}, all of the STR lengths \eqref{eqn:str_length} are zero, which implies degeneracy: $\mathcal{F}_{\boldsymbol{\eta}, f, \rho}^{\textnormal{STR}}[X] = \{\boldsymbol{\eta}\} \in \mathbb{L}_{\boldsymbol{\eta}}$.

\begin{flushright}
	$\blacksquare$
\end{flushright}
\subsection{Proof of \cref{prop:extreme_case_bounds}}
\label{app:proofs:prop:extreme_case_bounds}
Consider any bounded objective function $f: \mathbb{X} \times \Xi \rightarrow \mathbb{R}^M$ and strongly dominated reference vector $\boldsymbol{\eta} \in \cap_{\mathbf{x} \in \mathbb{X}} \cap_{\mathbf{y} \in f(\mathbf{x}, \Xi)} \mathbb{D}_{\precprec}[\{\mathbf{y}\}]$. Let $y_{\mathbf{x}, \boldsymbol{\xi}}^{(m)} = f^{(m)}(\mathbf{x}, \boldsymbol{\xi})$ denote the $m$-the objective value at the input $\mathbf{x} \in \mathbb{X}$ and uncertain parameter $\boldsymbol{\xi} \in \Xi$. Then for the upper bounds, we have that
\begin{equation*}
	\sup_{\boldsymbol{\xi} \in \Xi} \sup_{\mathbf{x} \in \mathbb{X}} \min_{m=1,\dots,M} z^{(m)}_{\boldsymbol{\eta}, \boldsymbol{\lambda}}(y_{\mathbf{x}, \boldsymbol{\xi}}^{(m)})
	=
	\sup_{\mathbf{x} \in \mathbb{X}} \sup_{\boldsymbol{\xi} \in \Xi} \min_{m=1,\dots,M} z^{(m)}_{\boldsymbol{\eta}, \boldsymbol{\lambda}}(y_{\mathbf{x}, \boldsymbol{\xi}}^{(m)})
	\leq
	\sup_{\mathbf{x} \in \mathbb{X}} \min_{m=1,\dots,M} z^{(m)}_{\boldsymbol{\eta}, \boldsymbol{\lambda}}(\sup_{\boldsymbol{\xi} \in \Xi} y_{\mathbf{x}, \boldsymbol{\xi}}^{(m)})
\end{equation*}
for any positive direction $\boldsymbol{\lambda} \in \mathcal{S}^{M-1}_+$. Equivalently, this means that $u_{\boldsymbol{\eta}, \boldsymbol{\lambda}, f}^{\textnormal{PS}} = u_{\boldsymbol{\eta}, \boldsymbol{\lambda}, f}^{\textnormal{STR}} \leq u_{\boldsymbol{\eta}, \boldsymbol{\lambda}, f}^{\textnormal{RTS}}$ for any positive direction $\boldsymbol{\lambda} \in \mathcal{S}^{M-1}_+$. Similarly, for the lower bounds, we have that
\begin{equation*}
	\sup_{\mathbf{x} \in \mathbb{X}} \min_{m=1,\dots,M} z^{(m)}_{\boldsymbol{\eta}, \boldsymbol{\lambda}}(\inf_{\boldsymbol{\xi} \in \Xi} y_{\mathbf{x}, \boldsymbol{\xi}}^{(m)})
	=
	\sup_{\mathbf{x} \in \mathbb{X}} \inf_{\boldsymbol{\xi} \in \Xi} \min_{m=1,\dots,M} z^{(m)}_{\boldsymbol{\eta}, \boldsymbol{\lambda}}(y_{\mathbf{x}, \boldsymbol{\xi}}^{(m)})
	\leq
	\inf_{\boldsymbol{\xi} \in \Xi} \sup_{\mathbf{x} \in \mathbb{X}} \min_{m=1,\dots,M} z^{(m)}_{\boldsymbol{\eta}, \boldsymbol{\lambda}}(y_{\mathbf{x}, \boldsymbol{\xi}}^{(m)})
\end{equation*}
for any positive direction $\boldsymbol{\lambda} \in \mathcal{S}^{M-1}_+$. Equivalently, this means that $- l_{\boldsymbol{\eta}, \boldsymbol{\lambda}, f}^{\textnormal{PS}} \leq - l_{\boldsymbol{\eta}, \boldsymbol{\lambda}, f}^{\textnormal{STR}} = - l_{\boldsymbol{\eta}, \boldsymbol{\lambda}, f}^{\textnormal{RTS}}$ for any positive direction $\boldsymbol{\lambda} \in \mathcal{S}^{M-1}_+$. The final result is then obtained by combining these two results together and using the fact that $l_{\boldsymbol{\eta}, \boldsymbol{\lambda}, f}^{\textnormal{PS}} \leq u_{\boldsymbol{\eta}, \boldsymbol{\lambda}, f}^{\textnormal{PS}}$ for any positive direction $\boldsymbol{\lambda} \in \mathcal{S}^{M-1}_+$.

\begin{flushright}
	$\blacksquare$
\end{flushright}
\subsection{Hypervolume of a Pareto front surface}
\label{app:proofs:hypervolume}
We now recall a useful result in \cref{lemma:hypervolume}, which states that the hypervolume of any Pareto front surface can be written as an integral over the set of positive unit vectors. A similar result was presented in the work by \cite{deng2019itec} and \cite{zhang2020icml} for the standard hypervolume indicator \eqref{eqn:hypervolume}.
\begin{lemma}
	[Hypervolume of a Pareto front surface] Consider a reference vector $\boldsymbol{\eta} \in \mathbb{R}^M$ and any Pareto front surface $A \in \mathbb{Y}_{\boldsymbol{\eta}}^*$, then the hypervolume of this polar surface can be written as
	\label{lemma:hypervolume}
	\begin{align*}
		\int_{\mathbb{R}^M} \mathbbm{1}[\mathbf{z} \in \mathbb{D}_{\preceq}[A] \cap \mathbb{D}_{\succsucc}[\{\boldsymbol{\eta}\}] ] d\mathbf{z}
		&= \mathbb{E}_{\boldsymbol{\lambda} \sim \textnormal{Uniform}(\mathcal{S}_+^{M-1})}
		[\tau^{\textnormal{HV}}(\ell_{\boldsymbol{\eta}, \boldsymbol{\lambda}}[A])].
	\end{align*}
\end{lemma}
\paragraph{Proof} By rewriting the volume integral in polar coordinates, we obtain the equation
\begin{align*}
	\int_{\mathbb{R}^M} \mathbbm{1}[\mathbf{z} \in \mathbb{D}_{\preceq}[A] \cap \mathbb{D}_{\succsucc}[\{\boldsymbol{\eta}\}]] d\mathbf{z}
	= c_M \int_{\boldsymbol{\lambda} \in \mathcal{S}^{M-1}_+} \ell_{\boldsymbol{\eta}, \boldsymbol{\lambda}}[A]^M d\boldsymbol{\phi}(\boldsymbol{\lambda}),
\end{align*}
where $\boldsymbol{\phi}$ denotes uniform measure over the positive unit vector $\boldsymbol{\lambda} \in \mathcal{S}^{M-1}_+$, whilst $c_M = \int_{\boldsymbol{\lambda} \in \mathcal{S}^{M-1}_+} d\boldsymbol{\phi}(\boldsymbol{\lambda}) = \pi^{M/2}/(2^M \Gamma(M/2 + 1))$ is the volume of the unit hypersphere lying in the positive orthant.
\begin{flushright}
	$\blacksquare$
\end{flushright}

\subsection{Proof of \cref{prop:rts_hypervolume} and \cref{prop:str_hypervolume}}
\label{app:proofs:prop:robust_hypervolume}
The proof of \cref{prop:rts_hypervolume} (or \cref{prop:str_hypervolume}) follows directly from the definition of the RTS front \eqref{eqn:rts_front} (or STR front \eqref{eqn:str_front}), together with \cref{prop:rts_front} (or \cref{prop:str_front}), \cref{lemma:hypervolume} and the monotonicity of the hypervolume transformation. 

\begin{flushright}
	$\blacksquare$
\end{flushright}
%
\bibliographystyle{plainnat}
\bibliography{ms}
\end{document}